\documentclass[preprint, 3p, 12pt]{elsarticle}

\usepackage{amssymb}

\usepackage[utf8]{inputenc}
\usepackage{textcomp}
\usepackage{ulem} 
\usepackage{graphicx}

\usepackage{amsmath}
\usepackage[version=4]{mhchem}
\usepackage{siunitx}
\usepackage{longtable,tabularx}
\usepackage{color}
\usepackage{mathtools}
\usepackage{float}
\usepackage{caption,subcaption}
\usepackage{tikz}
\usepackage[many]{tcolorbox}
\usepackage{color}
\usepackage{epstopdf}
\usepackage{algorithm}
\usepackage{algorithmic}
\usepackage{setspace}
\usepackage{mathrsfs}

\usepackage{hyperref}
\usepackage{ulem}

\def\1{\boldsymbol{1}}

\newtheorem{lemma}{Lemma}

\journal{Arxiv}

\begin{document}

\begin{frontmatter}


\title{Real-time optimal control for attitude-constrained solar sailcrafts via neural networks
}

\author{Kun Wang\fnref{label1}}
\ead{wongquinn@zju.edu.cn}
\author{Fangmin Lu\fnref{label1}}
\ead{lufangmin@zju.edu.cn}
\author{Zheng Chen\fnref{label1,label2} \corref{cor1}}
\ead{z\_chen@zju.edu.cn}
\author{Jun Li\fnref{label1,label2}}
\ead{lijun.uv@zju.edu.cn}


\cortext[cor1]{Corresponding author}
\affiliation[label1]{organization={School of Aeronautics and Astronautics, Zhejiang University}, 
            city={Hangzhou},
            postcode={310027},
            state={Zhejiang},
            country={P.R. China}}

\affiliation[label2]{organization={Huanjiang Lab}, 
            city={Zhuji},
            postcode={311816},
            state={Zhejiang},
            country={P.R. China}}





\begin{abstract}
This work is devoted to generating optimal guidance commands in real time for attitude-constrained solar sailcrafts in coplanar circular-to-circular interplanetary transfers. Firstly, a nonlinear optimal control problem is established, and necessary conditions for optimality are derived by the Pontryagin's Minimum Principle. Under 
some \textcolor{black}{mild} assumptions, the attitude constraints are rewritten as control constraints, which are replaced by a saturation function so that a parameterized system is formulated. \textcolor{black}{This allows one to generate an optimal trajectory via solving an initial value problem, making it 
efficient to collect a dataset containing optimal samples, which are essential for training Neural Networks (NNs) to achieve real-time implementation}. However, the optimal guidance command may suddenly change from one extreme to another, resulting in discontinuous jumps that generally impair the NN’s approximation performance. To address this issue, we use two co-states that the optimal guidance command depends on, to detect discontinuous jumps. A procedure for preprocessing these jumps is then established, thereby ensuring that the preprocessed guidance command remains smooth \textcolor{black}{everywhere}. Meanwhile, the sign of one co-state is found to be sufficient to revert the preprocessed guidance command back into the original optimal guidance command. Furthermore, three NNs are built and trained offline, and they cooperate together to precisely generate the optimal guidance command in real time. Finally, numerical simulations are presented to demonstrate the developments of the paper. 
\end{abstract}



\end{frontmatter}


\section{Introduction}
\label{intro}
Unlike conventional spacecrafts that use chemical or electric propellant to produce thrust during a flight mission, solar sailing exploits the \textcolor{black}{Solar Radiation Pressure} (SRP) to propel the solar sailcraft. Although the resulting propulsive force is smaller than that of chemical- or electric-based propulsion systems, the SRP offers an almost ``infinite'' energy source that can be used to orient the solar sailcraft \cite{tsu1959interplanetary}. This makes it a very promising technology, especially for long-duration interplanetary transfer missions \cite{fu2016solar,gong2019review}. 
While the study on solar sailing dates back to the 1920s, substantial progress in the development of solar sailing has not been achieved until recently. Following the first demonstrations of solar sailcraft in orbit by JAXA’s IKAROS \cite{mori2010first} and NASA’s NanoSail-D2 \cite{johnson2011nanosail} in 2010, the successes of LightSail 1 and 2 \cite{spencer2021lightsail} have sparked renewed interests in this technology. These achievements have led to the planning of numerous solar sailcraft-based missions, including Solar Cruiser \cite{spencer2021lightsail} and OKEANOS \cite{takao2021solar}.

Designing the transfer trajectory for solar sailcrafts is one of the most fundamental problems. As the acceleration provided by the SRP is extremely small, the transfer time of \textcolor{black}{the} solar sailcraft in space is usually too long. Thus, it is fundamentally important to find the  time-optimal control so that the solar sailcraft can be steered to its \textcolor{black}{targeted orbit} \textcolor{black}{within} minimum time. This is essentially equivalent to addressing a minimum-time control problem, which is a conventional Optimal Control Problem (OCP).  Numerical methods for OCPs include indirect and direct methods \cite{rao2009survey}. 
The indirect method, based on the calculus of variations or Pontryagin’s Minimum Principle (PMP), transforms the OCP into a Two-Point Boundary-Value Problem (TPBVP) according to the necessary conditions for optimality. The resulting TPBVP is then typically solved by Newton-like iterative methods. 
The indirect method has been applied to solve the time-optimal orbital transfer problem for solar \textcolor{black}{sailcrafts} in Refs. \cite{sauer1976optimum,wood1982comment,kim2005symmetries}.
Mengali and Quarta \cite{mengali2005optimal} employed the indirect method to solve the three-dimensional time-optimal transfer problem for a non-ideal solar sailcraft with optical and parametric force models. Sullo {\it et al.} \cite{sullo2017low} embedded \textcolor{black}{the continuation method} into the indirect method, so that the time-optimal transfer of \textcolor{black}{the} solar sailcraft was found from a previously designed low-thrust transfer trajectory. 
Recently, the solar sail primer vector theory, combined with the indirect method, was employed to design the flight trajectory that minimizes the solar angle over a transfer with \textcolor{black}{a fixed time of flight} \cite{oguri2023indirect}.
On the other hand, direct methods reformulate the OCP (usually via direct shooting or pseudospectral methods) as a nonlinear programming problem, which can be solved using interior point or sequential quadratic programming method \cite{rao2009survey}. \textcolor{black}{Solar} sailcraft transfers to hybrid pole and quasi pole sitters on Mars and Venus were studied using a direct pseudospectral algorithm in Ref. \cite{vergaaij2019solar}. In Ref. \cite{fitzgerald2021characterizing}, the Gaussian quadrature pseudospectral optimization was combined with a slerped control parameterization method, and the transfer time was drastically reduced. 
Furthermore, a comparison study on application of indirect methods and direct methods to \textcolor{black}{the} time-optimal transfer for solar sailcrafts was conducted in Ref. \cite{caruso2019comparison}. In addition to these conventional methods, heuristic methods have also been utilized in the literature to solve the time-optimal transfer problem; see, e.g., 
Ref. \cite{dachwald2005optimization}.

Although the techniques mentioned in the previous paragraph are widely used, they are usually time-consuming and need appropriate initial guesses of co-state or state vector. Therefore, these methods  are typically not suitable for onboard systems with limited computational resources. To overcome this issue, shape-based methods have been developed. The basic idea of the shape-based methods is to describe the geometric shape of trajectory using a set of mathematical expressions with some tunable parameters. 
The parameters are usually adjusted so that the mathematical expressions \textcolor{black}{match} the required boundary conditions. 
In this regard, Peloni {\it et al.} \cite{peloni2016solar} expressed the trajectory as a function of four parameters for a multiple near-Earth-asteroid mission. A Radau pseudospectral method was leveraged to solve the resulting problem in terms of these four parameters. Since then, 
different shaping functions, such as B\'{e}zier curve \cite{huo2019electric} and Fourier series \cite{taheri2012shape,caruso2021optimal}, \textcolor{black}{have been proposed for transfer problems}.
In order to further cope with constraints on the propulsive acceleration, some additional shape-based functions have been developed \cite{peloni2018automated,caruso2020shape}. It is worth mentioning that shape-based methods, albeit computationally efficient, often provide suboptimal solutions. Consequently, they are usually employed to provide initial guesses for direct methods.

According to the preceding review, existing approaches for solar sailcraft transfer suffer from computational burden, convergence issues, or solution suboptimality. In fact, solution optimality plays an 
important role in space exploration missions \cite{izzo2023optimality}, and it is demanding for onboard systems to generate real-time solutions. Therefore, it is worth exploiting more effective methods capable of providing optimal solutions in real time. Thanks to the emerging achievements of artificial intelligence in recent decades, machine learning techniques have become a viable approach.
Generally, there are two quite distinct machine learning based approaches. The first class, known as Reinforcement Learning (RL) involves training an agent to behave in a potentially changing environment via repeatedly interactions, \textcolor{black}{and the goal is to maximize a carefully designed reward function}. Owing to excellent generalization abilities of deep Neural Networks (NNs), deep RL algorithms have demonstrated great success in learning policies to guide spacecrafts in transfer missions; see, e.g., Refs \cite{zavoli2021reinforcement,federici2021autonomous,holt2021optimal}. \textcolor{black}{Although the trained agent is generalizable to non-nominal conditions and robust to uncertainties}, significant effort may be required to design an appropriate reward function, in order to derive a solution that is very close to optimal \cite{gaudet2020deep}.

In contrast to RL, Behavioral Cloning (BC) aims to clone the observed expert behavior through training an NN on optimal state-action pairs. Such pairs are usually obtained \textcolor{black}{by solving deterministic OCPs via indirect or direct methods}. Recently, BC has been widely used to generate optimal guidance commands in real time, or reduce the computational time \textcolor{black}{in aerospace applications}, such as spacecraft landing \cite{sanchez2018real,cheng2020real}, hypersonic vehicle reentry \cite{chai2019six,shi2021onboard}, missile interception \cite{wang2022nonlinear}, low-thrust transfer \cite{izzo2021real}, end-to-end human-Mars entry, powered-descent, and landing problem \cite{you2022onboard}, minimum-fuel geostationary  station-keeping \cite{zhang2023minimum}, time-optimal interplanetary  rendezvous \cite{izzo2023neural}, 
and solar sailcraft transfer \cite{cheng2018real}. 
Specifically, in Ref. \cite{cheng2018real}, 
some minimum-time orbital transfer problems with random initial conditions were solved by the indirect method to formulate the dataset for the state-guidance command pairs.  Then, the dataset was used to train NNs within a supervised learning framework. 
Finally, the trained NNs were able to generate the optimal guidance command in real time. 
However, constraints on the solar sailcraft attitude were not considered in that paper. In practice, factors such as power generation and thermal control generally reduce the admissible variation range of the cone angle, which is equal to the sail attitude for a perfectly flat ideal sail \cite{he2014time,caruso2020effects,oguri2022solar}. In addition, the optimal guidance command for the solar sailcraft transfer problem may be discontinuous, which generally degrades the approximation performance of NNs. In fact, approximating discontinuous controls via NNs, even with many hidden layers and/or neurons, can be quite challenging, as shown by the numerical simulations in Refs. \cite{li2019neural,george2022use,origer2023guidance}.  

As a continuation effort to the work in Ref. \cite{cheng2018real}, this paper considers \textcolor{black}{a scenario where} the solar sailcraft's attitude is constrained, and an NN-based method will be developed for real-time generation of \textcolor{black}{the} optimal guidance command. Firstly, the time-optimal problem for the solar sailcraft with constraints on \textcolor{black}{the} attitude is formulated. Then, necessary conditions for optimal trajectories are established by employing the PMP. These conditions allow using a saturation function to approximate the optimal guidance law. 
Unlike conventional optimization-based approaches, such as indirect methods and direct methods which often suffer from convergence issues \cite{chen2019nonlinear}, we formulate a parameterized system to  facilitate the generation of datasets for training NNs. In this method, an optimal trajectory can be generated by solving an initial value problem instead of a TPBVP. As a consequence, a large number of optimal trajectories can be readily obtained.

However, discontinuous jumps in the optimal guidance command may defer one from obtaining perfect \textcolor{black}{approximation}. 
To resolve this issue, one viable method is to use \textcolor{black}{NNs} to approximate smooth co-states instead, as was done in Ref. \cite{parrish2018low}. Unfortunately, using co-states to train NNs is not reliable as the relationship between co-states and \textcolor{black}{the} optimal guidance command is highly nonlinear. This may result in magnified propagation errors even for small errors from the predicted co-states \cite{rubinsztejn2020neural}. To avoid this difficulty, we propose to employ two co-states to preprocess the guidance command, smoothing out the discontinuous jumps.  After preprocessing, the guidance command can be reverted to its original optimal form by examining the sign of one specific co-state. To achieve real-time generation of the optimal guidance command, we employ a system comprising three NNs that predict the optimal time of flight, the preprocessed guidance command, and one co-state. These three NNs cooperate together to generate the optimal guidance command \textcolor{black}{in real time}. 

The remainder of the paper is structured as follows. The OCP is formulated in Section \ref{Problem}. Section \ref{Characterizations} presents the optimality conditions derived from the PMP, and a parameterized system for optimal solutions is established. In Section \ref{OptimalSteering}, procedures for generating the dataset and preprocessing discontinuous jumps are introduced, and the scheme for generating the optimal guidance command in real time is detailed. Section \ref{NumericalSimulations} provides some numerical examples to demonstrate the developments of the paper. This paper finally concludes by Section \ref{Colus}.
\section{Problem Formulation}
\label{Problem}
\subsection{Solar Sailcraft Dynamics}
\begin{figure}[!htp]
\centering
\includegraphics[width = 0.4\linewidth]{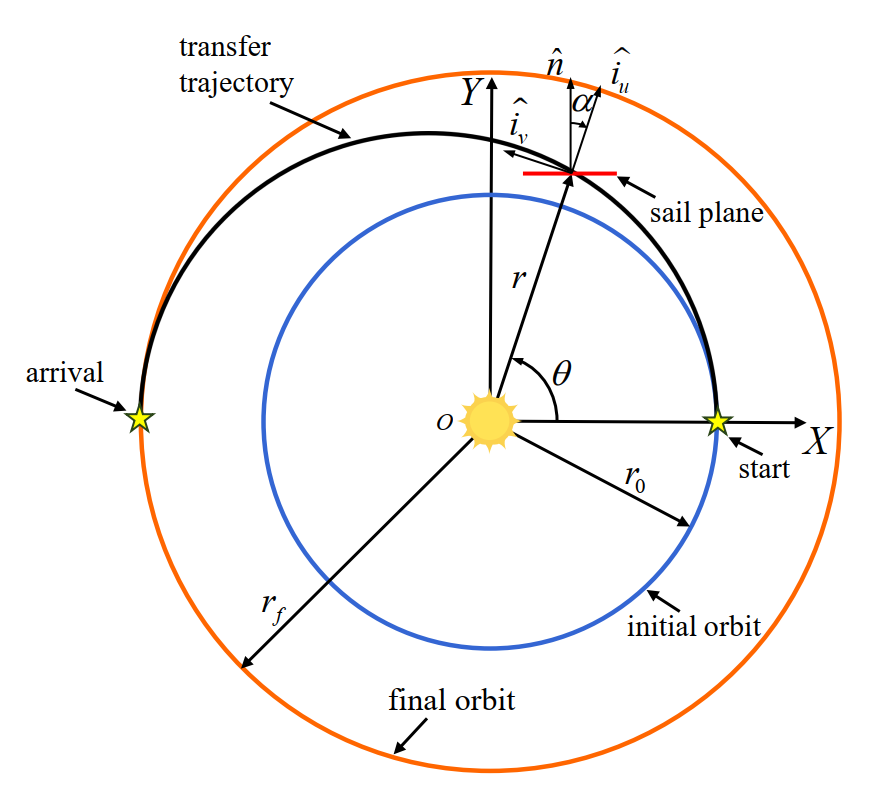}
\caption{Dynamics of an ideal solar sailcraft.}
\label{Fig:polar_reference}
\end{figure}

We consider a two-dimensional interplanetary transfer of an ideal solar sailcraft. Before proceeding, we present some units widely used in astrodynamics for better numerical conditioning. The distance is in units of Astronomical Unit (AU, \textcolor{black}{the} mean distance between the Earth and the Sun), and time is in units of Time Unit (TU, a time unit such that an object in a circular orbit at 1 AU would have a speed of 1 AU/TU). For clearer presentation, a heliocentric–ecliptic inertial frame $(O,X,Y)$ and a polar reference frame $(O,r,\theta)$ are used, as shown in Fig.~\ref{Fig:polar_reference}. The origin $O$ is located at the Sun's center; the $Y$ axis points towards the opposite direction of the vernal equinox, and the $X$ axis is specified by rotating the $Y$ axis 90 degrees clockwise in the ecliptic plane. $r \in \mathbb{R}_0^+$ denotes the distance from the origin to the solar sailcraft, and $\theta \in [0,2\pi]$ is the rotation angle of the Sun-solar sailcraft line measured counterclockwise in the polar reference frame from the positive $X$ axis. 

In the context of preliminary mission design, it is reasonable to make some assumptions. The Sun and solar sailcraft are treated as point masses.
The motion of the solar sailcraft is only influenced by the Sun’s gravitational attraction and the propulsive acceleration due to the SRP acting on the solar sailcraft. Other disturbing factors, such as the solar wind and light aberration are neglected. Then, the dimensionless dynamics of the solar sailcraft is given by \cite{kim2005symmetries}
\begin{align}
\begin{cases}
\dot{r}(t) = u(t),\\
\dot{\theta}(t) = \frac{v(t)}{r(t)},\\
\dot{u}(t) = \frac{\beta \cos^3\alpha(t)}{r^2(t)} + \frac{v^2(t)}{r(t)} - \frac{1}{r^2(t)},\\
\dot{v}(t) = \frac{\beta \sin\alpha(t) \cos^2\alpha(t)}{r^2(t)} - \frac{u(t)v(t)}{r(t)},
\label{EQ:dyna_equation}
\end{cases}
\end{align}
where $t \geq 0$ is the time; $u \in \mathbb{R}$ and $v \in \mathbb{R}$ is the radial and transverse speed in units of AU/TU, respectively; $\beta$ is the normalized characteristic acceleration; $\alpha$ represents the cone angle, which is the angle between the normal direction of sail plane $\hat{n}$ and the direction of incoming rays $\hat{i}_u$, measured positive in the clockwise direction from the normal direction of sail plane $\hat{n}$.
\subsection{Attitude Constraints}
Note that the radial component $a_u$ and transverse component $a_v$ of the ideal propulsive acceleration vector acting on the solar sailcraft are given by \cite{caruso2020effects}
\begin{align}
\begin{cases}
a_u(t) := \frac{\beta \cos^3\alpha(t)}{r^2(t)}, \\
a_v(t) := \frac{\beta \sin\alpha(t) \cos^2\alpha(t)}{r^2(t)}. 
\label{EQ:acc_equation}
\end{cases}
\end{align}
Based on the assumption that the ideal sail is approximated with a flat and rigid reflecting surface, the attitude \textcolor{black}{constraints} can be rewritten in terms of \textcolor{black}{the} cone angle as
\begin{align}
\alpha \in [-\phi_{max},\phi_{max}], 
\label{EQ:att_constraints}
\end{align}
in which $\phi_{max} \in (0,\frac{\pi}{2})$ is a given parameter depending on the minimum acceptable level of the electric power generated by the solar sailcraft \cite{caruso2020effects}. It essentially sets the limits for the solar sailcraft's orientation with respect to the incoming solar rays.

\begin{figure}[!htp]
\centering
\includegraphics[width = 0.8\linewidth]{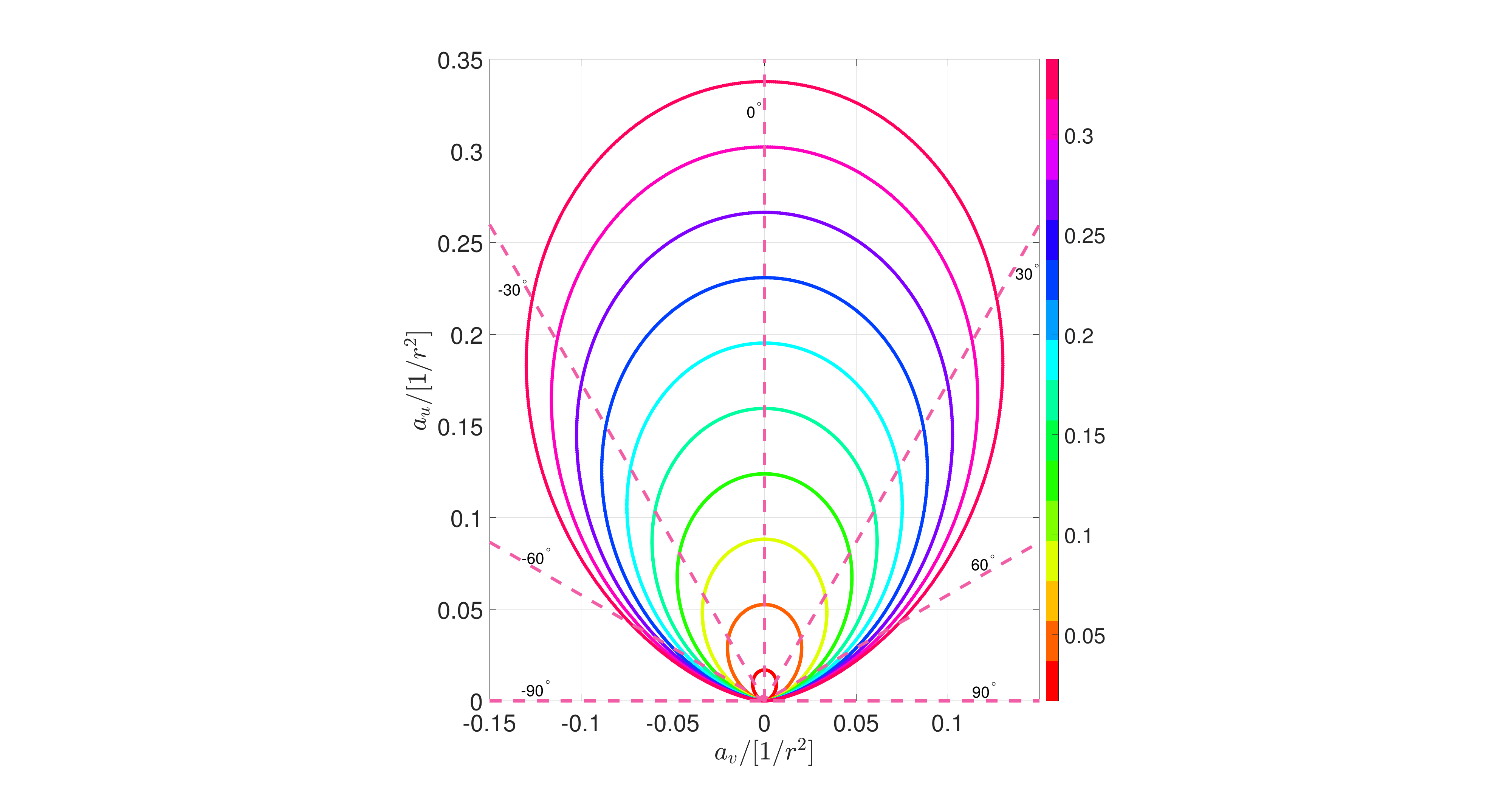}
\caption{Shape of the ideal sail force bubble.}
\label{Fig:force_bubble}
\end{figure}

The ``force bubble'' of the ideal solar mentioned in Ref. \cite{mengali2007refined} due to attitude constraints is displayed in Fig.~\ref{Fig:force_bubble}. The colormap represents the value of $\beta$ in $[0.01892, 0.3784]$ (corresponding to a characteristic acceleration domain of $\rm {[0.1, 2]~mm/s^2}$). Then, for a given $r$, it is clear that $\beta$ defines the size of the force bubble, whereas $\phi_{max}$ determines its actual shape, as shown by the pink-dashed lines. Specifically, 
\textcolor{black}{the propulsive acceleration is constrained to the radial direction if $\phi_{max} = 0$.}
\subsection{Formulation of the OCP}
Without loss of generality, consider an initial condition for the solar sailcraft at the initial time $t_0 = 0$ given by
\begin{align}
r(0) = r_0, \theta(0) = \theta_0, u(0) = u_0, v(0) = v_0.
\label{EQ:initial_condition}
\end{align}
The mission purpose is, through orientating the cone angle $\alpha$ subject to constraints in Eq.~(\ref{EQ:att_constraints}), to drive the solar sailcraft governed by Eq.~(\ref{EQ:dyna_equation}), into a coplanar circular orbit of radius $r_f$ with the terminal condition given by
\begin{align}
r(t_f) = r_f, u(t_f) = 0, v(t_f) = \frac{1}{\sqrt{r_f}}.
\label{EQ:final_condition}
\end{align}
The functional cost $J$ is to minimize the \textcolor{black}{arrival time (final time)} $t_f$, i.e.,
\begin{align}
J = \int_{0}^{t_f}1~dt.
\label{EQ:cost_function}
\end{align}  
\section{Parameterization of Optimal Trajectories}
\label{Characterizations}
In this section, we first derive the necessary conditions for optimal trajectories. Then, we formulate a parameterized system that enables the generation of an optimal trajectory via solving an initial value problem. For simplicity of presentation, important results and claims are written in lemmas with their proofs postponed to the appendix.
\subsection{Optimality Conditions}
Denote by $\boldsymbol{\lambda} = [\lambda_r, \lambda_\theta, \lambda_{u},\lambda_{v}]$ the co-state vector of the state vector $\boldsymbol{x} = [r, \theta, u,v]$. Then, the Hamiltonian for the OCP is expressed as
\begin{align}
\mathscr H = 1 + \lambda_r u + \lambda_\theta \frac{v}{r} + \lambda_{u}(\frac{\beta \cos^3\alpha}{r^2} + \frac{v^2}{r} - \frac{1}{r^2}) + \lambda_{v}
(\frac{\beta \sin\alpha \cos^2\alpha}{r^2} - \frac{u v}{r}).
\label{EQ:Ham_function}
\end{align}
According to the PMP \cite{Pontryagin}, we have
\begin{align}
\begin{cases}
\dot{\lambda}_r(t) = \lambda_\theta(t)\frac{v(t)}{r^2(t)}+\lambda_{u}(t)[2\beta\frac{\cos^3\alpha(t)}{r^3(t)}+\frac{{v}^2(t)}{r^2(t)}-\frac{2}{r^3(t)}]+\lambda_{v}(t)[2\beta \cos^2 \alpha(t) \frac{\sin \alpha(t)}{r^3(t)}-\frac{u(t)v(t)}{r^2(t)}],\\
\dot{\lambda}_\theta(t) = 0,\\
\dot{\lambda}_{u}(t) = -\lambda_r(t) + \frac{\lambda_{v}(t) v(t)}{r(t)},\\
\dot{\lambda}_{v}(t) = -\frac{\lambda_\theta(t)}{r(t)} - 2\frac{\lambda_{u}(t)v(t)}{r(t)} + \frac{\lambda_{v}(t)u(t)}{r(t)}.
\label{EQ:costate_function}
\end{cases}
\end{align}
Because $t_f$ is free, it holds
\begin{align}
\mathscr H(t) \equiv 0,  ~\forall~t \in [0,t_f].
\label{EQ:tf_law}
\end{align}
In addition, the terminal rotation angle $\theta(t_f)$ is not constrained, leading to
\begin{align}
\lambda_\theta(t_f) = 0.
\label{EQ:theta_law}
\end{align}
By the following lemma, we present the optimal guidance law.
\begin{lemma}\label{LE:optimal_control_law}
The optimal guidance law that minimizes $\mathscr H$ in Eq.~(\ref{EQ:Ham_function}) is
\begin{align}
\alpha^*=\text{Median} {[-\phi_{max},\bar{\alpha},\phi_{max}]}, \text{where}~\bar{\alpha} = \arctan~\frac{-3{\lambda_{u}}-\sqrt{9\lambda^2_{u}+8\lambda^2_{v}}}{4\lambda_{v}}.
\label{EQ:optimal_control}
\end{align}
\end{lemma}
The proof of this lemma is postponed to \ref{Appendix:A}.  

Notice that the optimal guidance law in Eq.~(\ref{EQ:optimal_control}) may not be smooth at some isolated points. In fact, the optimal guidance law in Eq.~(\ref{EQ:optimal_control}) can be approximated in a compact form as ~\cite{avvakumov2000boundary} 
\begin{align}
\alpha^* \approx \alpha^*(\delta)= \frac{1}{2}~[\sqrt{\left(\bar{\alpha}+\phi_{max}\right)^2+\delta}-\sqrt{\left(\bar{\alpha}-\phi_{max}\right)^2+\delta}].
\label{EQ:optimal_control_smooth}
\end{align}
where $\delta$ is a small non-negative number. It is clear that Eq.~(\ref{EQ:optimal_control_smooth}) is equivalent to Eq.~(\ref{EQ:optimal_control}) if $\delta=0$. If $\delta>0$ is sufficiently small, Eq.~(\ref{EQ:optimal_control_smooth}) acts like a smoothing filter function. 

Note that the initial rotation angle $\theta_0$ has no effect on the solutions because the transfer trajectory is rotatable \cite{cheng2018real}. For brevity, a triple $(r(t),u(t),v(t))$ for $t \in [0,t_f]$ is said to be an optimal trajectory if all the necessary conditions in Eqs.~(\ref{EQ:costate_function}), (\ref{EQ:tf_law}), (\ref{EQ:theta_law}), and (\ref{EQ:optimal_control}) are met. In order to generate the optimal guidance command in real time via NNs, a training dataset containing a large number of optimal trajectories is required. In this regard, one viable approach is to employ some root-finding algorithms to solve the TPBVP
\begin{align}
[r(t_f)-r_f,u(t_f),v(t_f)-\frac{1}{\sqrt{r_f}},\mathscr{H}(t_f)]=\boldsymbol{0}.
\label{EQ:TPBVP_law}
\end{align}
Nevertheless, this procedure is usually time-consuming and suffers from convergence issues. In the next subsection, we present a parameterized system so that an optimal trajectory can be readily obtained by solving an initial value problem \textcolor{black}{instead}.
\subsection{Parameterized System}
Define a new independent variable $\tau$ as below
\begin{align}
\tau = t_f - t, t \in [0,t_f].
\end{align}
Let us establish a first-order ordinary differential system 
\begin{align}
\begin{cases}
\dot{R}(\tau) = -U(\tau),\\
\dot{U}(\tau) = -\frac{\beta \cos^3\hat{\alpha}(\tau)}{R^2(\tau)} - \frac{V^2(\tau)}{R(\tau)} + \frac{1}{R^2(\tau)},\\
\dot{V}(\tau) = -\frac{\beta \sin\hat{\alpha}(\tau) \cos^2\hat{\alpha}(\tau)}{R^2(\tau)} + \frac{U(\tau)V(\tau)}{R(\tau)},\\
\dot{\lambda}_R(\tau) = -\lambda_{U}(\tau)[2\beta\frac{\cos^3\hat{\alpha}(\tau)}{R^3(\tau)}\frac{{V}^2(\tau)}{R^2(\tau)}-\frac{2}{R^3(\tau)}]-\lambda_{V}(\tau)[2\beta \cos^2 \hat{\alpha}(\tau) \frac{\sin \hat{\alpha}(\tau)}{R^3(\tau)}-\frac{U(\tau)V(\tau)}{R^2(\tau)}],\\
\dot{\lambda}_{U}(\tau) = \lambda_R(\tau) - \frac{\lambda_{V}(\tau) V(\tau)}{R(\tau)},\\
\dot{\lambda}_{V}(\tau) =  2\frac{\lambda_{U}(\tau)V(\tau)}{R(\tau)} - \frac{\lambda_{V}(\tau)U(\tau)}{R(\tau)},
\label{EQ:new_system}
\end{cases}
\end{align}
where $(R, U, V, \lambda_R, \lambda_{U}, \lambda_{V}) \in  \mathbb{R}_0^+ \times \mathbb{R}^5$, and $\hat{\alpha}$ is defined as
\begin{align}
\hat{\alpha}=\frac{1}{2}~[\sqrt{\left(\bar{\alpha}+\phi_{max}\right)^2+\delta}-\sqrt{\left(\bar{\alpha}-\phi_{max}\right)^2+\delta}]~\rm{with}~
 \bar{\alpha} = \arctan~{\frac{-3 \lambda_{U}-\sqrt{9\lambda^2_{U}+8\lambda^2_{V}}}{4\lambda_{V}}}. 
\label{EQ:new_control}
\end{align}

The initial value at $\tau = 0$ for the system in Eq.~(\ref{EQ:new_system}) is set as
\begin{align}
R(0) = R_0, U(0) = 0, V(0) = \frac{1}{\sqrt{R_0}}, \lambda_{R}(0) = \lambda_{R_0}, \lambda_{U}(0) = \lambda_{U_0}, \lambda_{V}(0) = \lambda_{V_0}.
\label{para_initial}
\end{align}
The value of $\lambda_{V_0}$ satisfies  \textcolor{black}{the following equation}
\begin{equation}
        \begin{split}
1 + \lambda_{R}(0) U(0) +  \lambda_{U}(0)(\frac{\beta \cos^3\hat{\alpha}(0)}{R^2(0)} + \frac{V^2(0)}{R(0)} - \frac{1}{R^2(0)}) \\ 
+ \lambda_{V}(0)(\frac{\beta \sin\hat{\alpha}(0) \cos^2\hat{\alpha}(0)}{R^2(0)} - \frac{U(0) V(0)}{R(0)})= 0.
        \end{split}
\label{EQ:solve_equation}
\end{equation}
Substituting Eq.~(\ref{para_initial}) into Eq.~(\ref{EQ:solve_equation}) leads to
\begin{align}
\textcolor{black}{
1 +  \lambda_{U_0}\frac{\beta \cos^3\hat{\alpha}(0)}{R^2_0} + \lambda_{V_0}
\frac{\beta \sin\hat{\alpha}(0) \cos^2\hat{\alpha}(0)}{R^2_0} = 0.}
\label{EQ:solve_equation_simple}
\end{align}
\textcolor{black}{
In view of Eqs.~(\ref{EQ:new_control}) and (\ref{EQ:solve_equation_simple}), it is clear that for a given $\lambda_{U_0}$, the value for $\lambda_{V_0}$ can be numerically determined by solving Eq.~(\ref{EQ:solve_equation_simple}).}
By the following lemma, we shall show that an optimal trajectory can be generated by solving an initial value problem based on
the parameterized system in Eq.~(\ref{EQ:new_system}).

\begin{lemma}\label{LE:optimal_trajectory_lemma}
\textcolor{black}{For any given constants $\lambda_{R_0}$ and $\lambda_{U_0}$}, a fixed \textcolor{black}{final} time $t_f$, and \textcolor{black}{a parameter} $\lambda_{V_0}$ \textcolor{black}{to be} determined by \textcolor{black}{solving} Eq.~(\ref{EQ:solve_equation_simple}), denote by 
\begin{align*}
\mathcal{F}:=
(R(\tau,\lambda_{R_0},\lambda_{U_0}),U(\tau,\lambda_{R_0},\lambda_{U_0}),V(\tau,\lambda_{R_0},\lambda_{U_0}),\\
\lambda_R(\tau,\lambda_{R_0},\lambda_{U_0}),\lambda_{U}(\tau,\lambda_{R_0},\lambda_{U_0}),
\lambda_{V}(\tau,\lambda_{R_0},\lambda_{U_0})) \in \mathbb{R}^+ \times \mathbb{R}^5.
\end{align*}
the solution of the $(\lambda_{R_0},\lambda_{U_0})$-parameterized system in
Eq.~(\ref{EQ:new_system}) with the initial value specified in Eq.~(\ref{para_initial}). Define $\mathcal{F}_p$ as
\begin{align*}
\mathcal{F}_p:=
\{R(\tau,\lambda_{R_0},\lambda_{U_0}),U(\tau,\lambda_{R_0},\lambda_{U_0}),V(\tau,\lambda_{R_0},\lambda_{U_0}),\tau |\\
 (R(\tau,\lambda_{R_0},\lambda_{U_0}),U(\tau,\lambda_{R_0},\lambda_{U_0}),V(\tau,\lambda_{R_0},\lambda_{U_0}),\tau) \in \mathcal{F}\}.
\end{align*}
Then $\mathcal{F}_p$ represents the solution space of an optimal trajectory starting from a circular orbit with a radius of $R_0$.
\end{lemma}

The proof of this lemma is postponed to \ref{Appendix:A}.

With the parameterized system established, obtaining an optimal trajectory becomes a straightforward process that involves solving an initial value problem, rather than tackling a TPBVP.
Once the optimal trajectory $\mathcal{F}_p$ is obtained, the corresponding optimal guidance command $\hat{\alpha}$ can be determined through the definition of $\mathcal{F}$ and reference to Eq.~(\ref{EQ:new_control}). 
To maintain clarity and prevent any ambiguity in notation,  we continue to employ $(r,u,v)$ to represent the flight state and \textcolor{black}{we use $\tau$} to represent the optimal time of flight along the optimal trajectory $\mathcal{F}_p$. Since the optimal time of flight is of great importance for preliminary mission design, evaluating the optimal time of flight without the need for the exact guidance law is quite attractive, as was done in Refs. \cite{cheng2020real,izzo2021real}. In this paper, 
we include \textcolor{black}{$\tau$} into the flight state, and it can be evaluated via a trained NN, as shown by our scheme in Subsection \ref{SchemeGuidanceCommands}.
 Likewise, the optimal guidance command is again denoted by $\alpha$.
Then we shall define $f$ as the nonlinear mapping from the flight state \textcolor{black}{$(r,u,v,\tau)$} to the optimal guidance command $\alpha$, i.e.,
\begin{align*}
\textcolor{black}{
f:=(r,u,v,\tau) \to \alpha}.
\end{align*}

According to the universal approximation theorem \cite{hornik1989multilayer}, employing a large number of sampled data that map the flight state \textcolor{black}{$(r,u,v,\tau)$} to the optimal guidance command $\alpha$ to train feedforward NNs enables the trained networks to accurately approximate the nonlinear mapping $f$. 
Before generating the dataset, a \it{nominal trajectory} \rm{is typically required}. Then, by applying perturbations over $\lambda_{R_0}$ and $\lambda_{U_0}$ on the  \it{nominal trajectory}, \rm{lots of optimal trajectories containing optimal state-guidance command pairs can be aggregated.} With this dataset, the optimal guidance command can be generated in real time within the supervised learning framework, as we will detail in the next section. 
\section{Generation of Optimal Guidance Commands in Real Time via NNs}\label{OptimalSteering}
In this section, we begin by outlining the procedure for generating the training dataset. As our analysis progresses, we observe a notable characteristic of the optimal guidance command, \textcolor{black}{that is,} it may abruptly change from one extreme to another, resulting in discontinuous jumps. This behavior poses a challenge, as the approximation performance of the NN typically deteriorates under such circumstances, a phenomenon seen in Refs. \cite{li2019neural,george2022use,origer2023guidance}.
To address this issue, we introduce a preprocessing method that employs two co-states to detect the discontinuous jump. With the saturation function defined in Eq.~(\ref{EQ:optimal_control_smooth}), the preprocessed guidance command becomes smooth everywhere. 
Following this procedure, we present a scheme for generating optimal guidance commands in real time via NNs. 
\subsection{Generation of the Training Dataset}
\subsubsection{Nominal Trajectory}
\label{nominal}
Consider a solar sailcraft subject to attitude constraints $\alpha \in [-\frac{\pi}{3}, \frac{\pi}{3}]$, starting from Earth's orbit around the Sun, which is a circular orbit with radius of 1 AU (1 AU = $1.4959965\times 10^{11}$ m), i.e.,
\begin{align}
r(0) = 1, \theta(0) = 0, u(0) = 0, v(0) = 1,
\label{EQ:initial_condition_nominal}
\end{align}
and the mission is to fly into Mars' orbit around the Sun, which is 
a circular orbit with radius of 1.524 AU, i.e., 
\begin{align}
r(t_f) = 1.524, u(t_f) = 0, v(t_f) = \frac{1}{\sqrt {1.524}}.
\label{EQ:final_condition_nominal}
\end{align}

The normalized
characteristic acceleration $\beta$ is set as a constant of 0.16892, corresponding to a characteristic acceleration of 1 $\rm{mm/s^2}$.
To demonstrate the effects of the saturation function in Eq.~(\ref{EQ:optimal_control_smooth}), the optimal guidance laws in Eqs.~(\ref{EQ:optimal_control}) and (\ref{EQ:optimal_control_smooth}) are embedded into the shooting function defined by Eq.~(\ref{EQ:TPBVP_law}). To enhance the convergence of the indirect method, we initially consider cases without attitude constraints and subsequently apply a homotopy technique to accommodate attitude constraints. 

\begin{figure}[!htp]
\centering
\includegraphics[width = 0.7\linewidth]{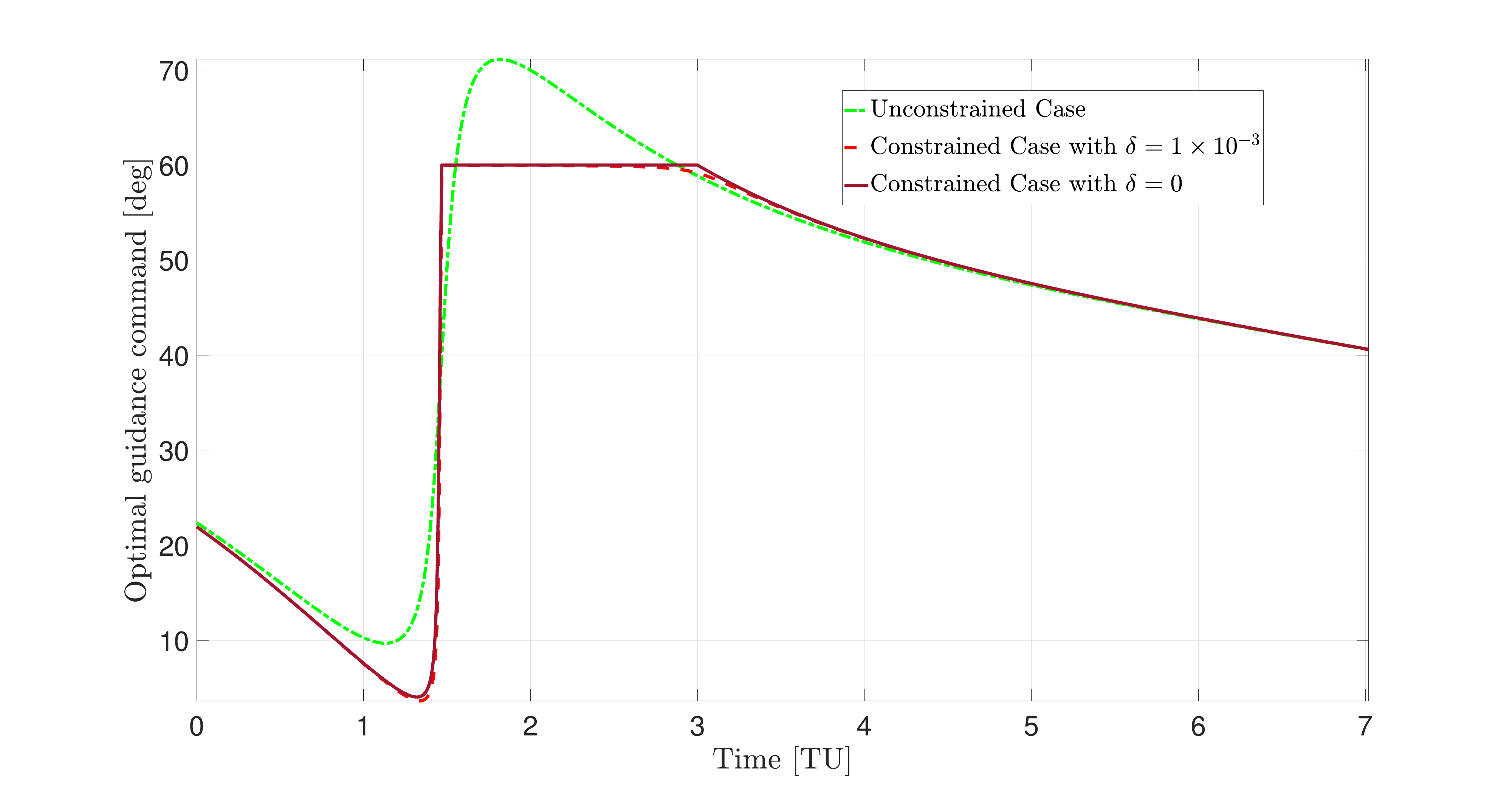}
\caption{Optimal guidance command profiles for the unconstrained and constrained cases.}
\label{Fig:norm_angle_two}
\end{figure}

Fig.~\ref{Fig:norm_angle_two}
shows the optimal guidance command profiles for the unconstrained and constrained cases with two smoothing constants. It is evident that the constrained cases contain rapid changes in the guidance command. Additionally, the guidance command with $\delta = 1 \times 10^{-3}$ displays smoother behavior compared to the case with $\delta = 0$. Besides, the arrival time for the unconstrained case is 7.01204 TU ( 1 TU $=5.0226757\times 10^{6}$ s), whereas this time span extends to $7.01818$ TU for the constrained case with $\delta = 0$. Regarding $\delta = 1 \times 10^{-3}$, the arrival time is $7.01838$ TU, indicating that a smoothing constant of $\delta = 1 \times 10^{-3}$ leads to a minor functional penalty of merely $2 \times 10^{-4}$ TU. Consequently, we designate the transfer trajectory with $\delta=1 \times 10^{-3}$ as the \it{nominal trajectory} \rm{henceforth}.
\subsubsection{Dataset}
Denote by $\lambda^*_{r_f}$, $\lambda^*_{{u}_f}$ the final values along the \it{nominal trajectory} \rm{for the co-state $\lambda_{r}$ and $\lambda_{u}$, respectively.} \textcolor{black}{They are computed as  $\lambda^*_{r_f}=-21.51$ and $\lambda^*_{{u}_f} = 8.54$.}
Consider a new set of values given by
\begin{align}
\lambda'_{R_0} = \lambda^*_{r_f} + \delta \lambda_R, \lambda'_{U_0} = \lambda^*_{{u}_f} + \delta \lambda_{U},
\label{EQ:newset}
\end{align}
\textcolor{black}{where the perturbations $\delta \lambda_R$ and $\delta \lambda_{U}$ are chosen such that the resulting $\lambda'_{R_0}$ and $\lambda'_{U_0}$ are uniformly distributed in the interval $[-23,-20]$ and $[5,11]$, respectively}. Set $ \lambda_{R_0} = \lambda'_{R_0}$,$ \lambda_{U_0} = \lambda'_{U_0}$, and calculate the value of $\lambda_{V_0}$ by solving Eq.~(\ref{EQ:solve_equation_simple}).  Then, $\mathcal F$ can be obtained by propagating the $(\lambda_{R_0},\lambda_{U_0})$-parameterized system with the initial value specified in Eq.~(\ref{para_initial}). Define an empty set $\mathcal D$, and insert $\mathcal {F}_p$ depending on the pair $(\lambda_{R_0},\lambda_{U_0})$, into $\mathcal D$ until the perturbation process described in Eq.~(\ref{EQ:newset}) ends. In this context, we can amass a dataset $\mathcal D$ comprising the essential optimal state-guidance command pairs required for NN training. 

Because we are specifically dealing with orbit transfer missions from Earth to Mars in this study, any optimal trajectory with $R(\tau) > 1.54$ AU for $\tau \in [0,t_f]$ ($t_f$ is fixed as 10) is classified as beyond the region of interest and excluded from $\mathcal D$. Moreover, to explore the NN's generalization ability, the \it{nominal trajectory} \rm{is excluded from $\mathcal D$}, as was done in Ref. \cite{izzo2021real}. Ultimately, a total of $18,324$ optimal trajectories are acquired, and they are evenly sampled, resulting in $6,416,613$ pairs. 
\subsection{Preprocessing Discontinuous Jumps}\label{Jumps}
\begin{figure}[!htp]
\centering
\includegraphics[width = 0.7\linewidth]{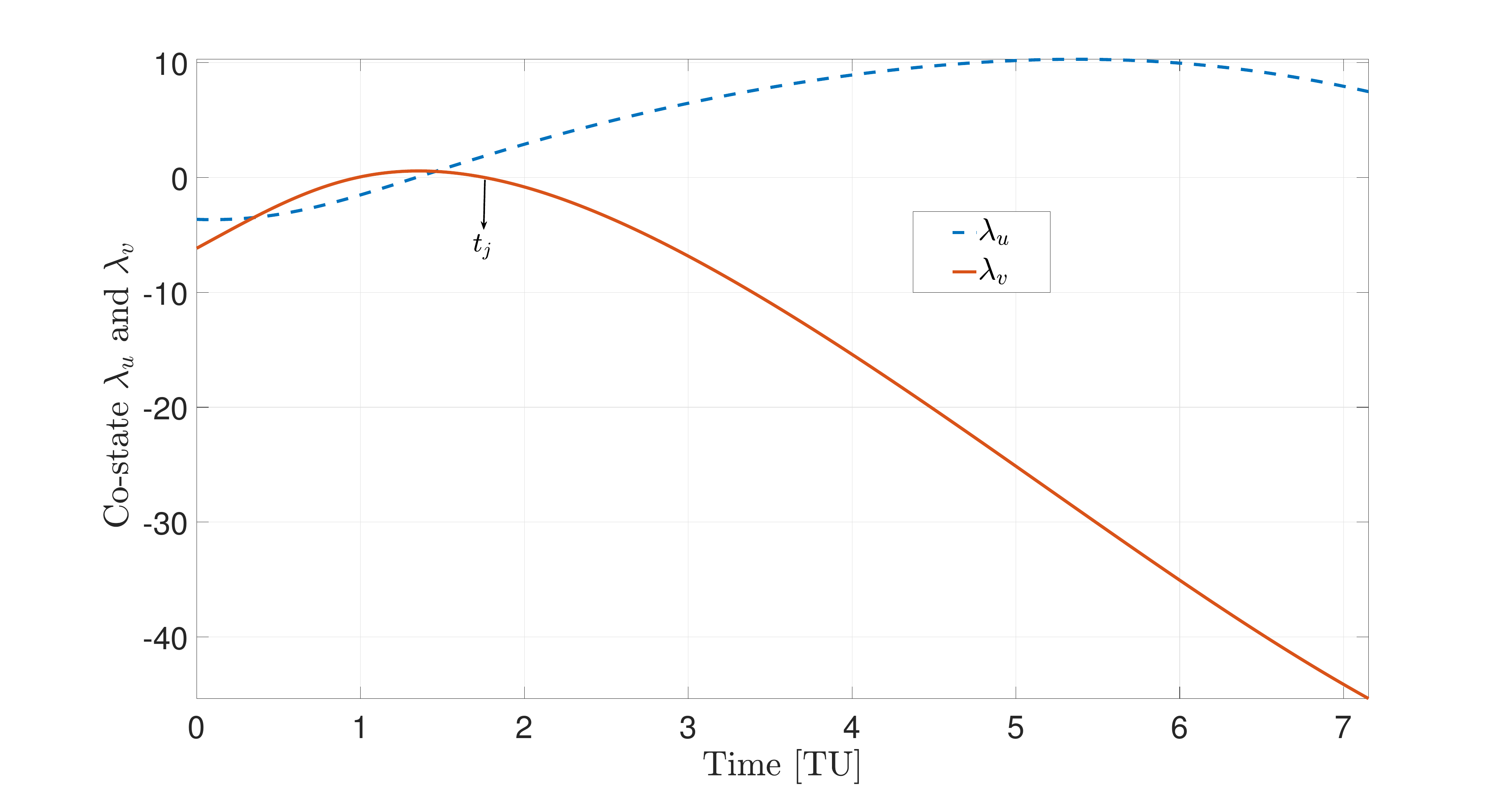}
\caption{Two co-state profiles along an optimal trajectory.}
\label{Fig:costate_ex}
\end{figure}

While generating the dataset, we encounter instances where the optimal guidance command displays discontinuous jumps. Remember that these jumps can be detected using $\lambda_{u}$ and $\lambda_{v}$, as elaborated in \ref{Appendix:A}. For a representative optimal trajectory defined in $[0,t_f]$, spotting a discontinuous jump point $t_j$ is straightforward, as illustrated in Fig.~\ref{Fig:costate_ex}. 
The optimal trajectory is then split into two segments, \textcolor{black}{i.e.,} $[0,t_j]$ and $[t_j,t_f]$. Assume that there is only one jump in this scenario, we can derive a preprocessed smooth guidance command, denoted by $\alpha_{pre}$, using the following approach
\begin{align}
   \alpha_{pre} = 
   \begin{cases}
      -\alpha, \text{for}~ {t \in [0, t_j]}\\
      ~~\alpha, \text{for}~ {t \in [t_j, t_f]}
   \end{cases}
   \label{regulaeise}
\end{align}

\begin{figure}[!htp]
\centering
\includegraphics[width = 0.7\linewidth]{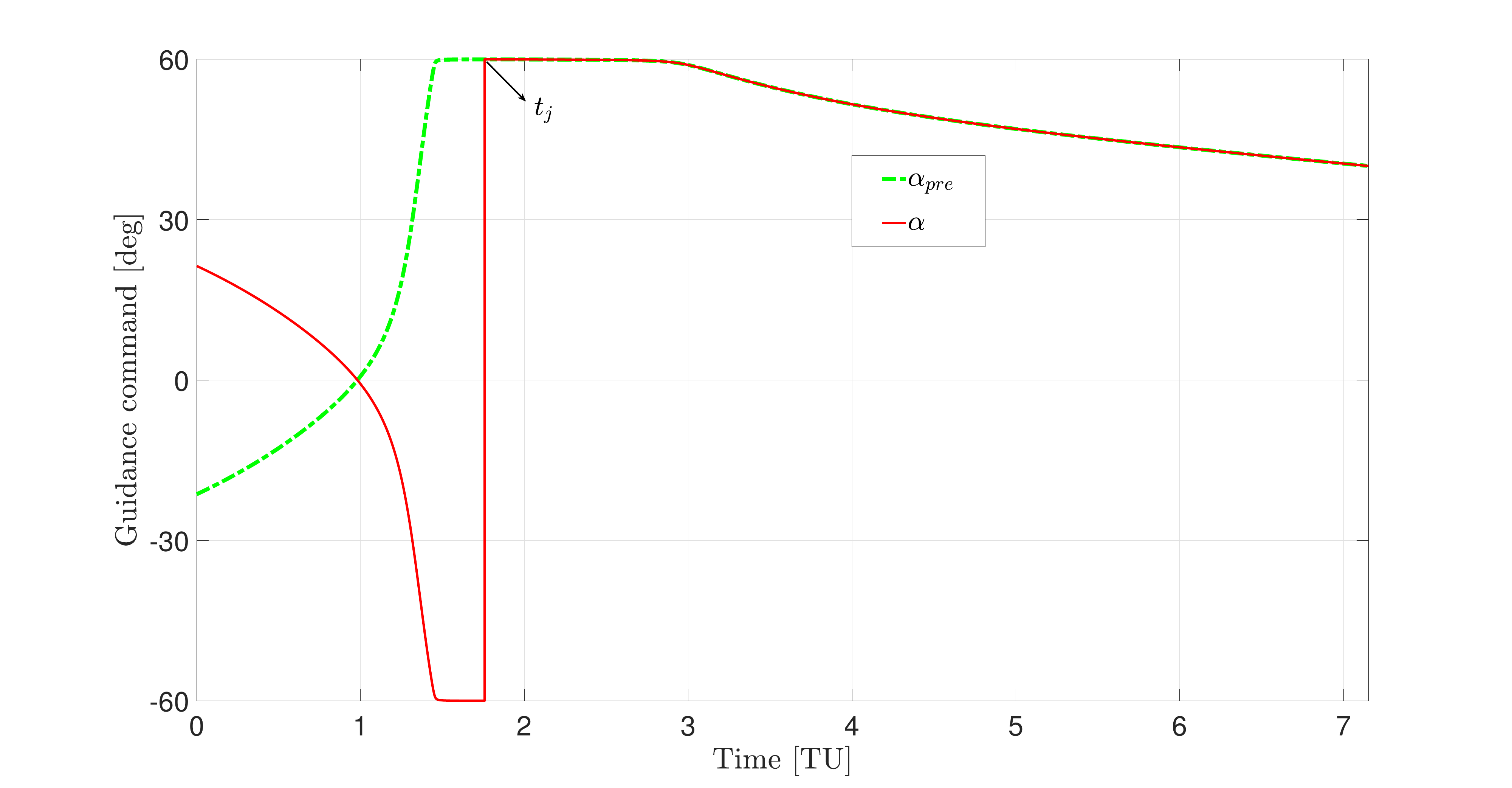}
\caption{Preprocessed and original optimal guidance command profiles along an optimal trajectory.}
\label{Fig:alpha_ex}
\end{figure}

A representative example is provided in Fig.~\ref{Fig:alpha_ex}. 
It is evident that, thanks to the saturation function in Eq.~(\ref{EQ:optimal_control_smooth}) and the preprocessing procedure outlined in Eq.~(\ref{regulaeise}), $\alpha_{pre}$ now demonstrates smooth behavior along the entire optimal trajectory. Conversely, if the discontinuous jump does not appear along the optimal trajectory, it is clear that $\alpha_{pre} = \alpha$ \textcolor{black}{always holds} for $t \in [0,t_f]$. Moreover, we observe that multiple discontinuous jumps are rare and primarily occur in orbital transfers with unusually \textcolor{black}{long} time of flight. Hence, in this study, we limit our consideration to orbital transfers with a maximum of one discontinuous jump.

With the parameterized system, a sample of dataset $\mathcal{D}_{\alpha}:=\left\{(r,u,v,\tau),\alpha \right\}$ can be extracted from the dataset $\mathcal{D}$. Thanks to the preprocessing procedure aforementioned, we can transform the dataset $\mathcal{D}_{\alpha}$, which could contain discontinuous jumps, into a refined dataset $\mathcal{D}_{\alpha_{pre}}:=\left\{(r,u,v,\tau),\alpha_{pre} \right\}$ \textcolor{black}{exclusively} consisting of smooth guidance commands. Consequently, a well-trained NN $\mathcal{N}_{\alpha_{pre}}$ based on $\mathcal{D}_{\alpha_{pre}}$ is anticipated to yield reduced \textcolor{black}{approximation} errors compared to an NN $\mathcal{N}_\alpha$ trained on $\mathcal{D}_{\alpha}$ \cite{izzo2023optimality}.

For the sake of completeness, we will elucidate the process of reverting the output $\alpha_{pre}$ from $\mathcal{N}_{\alpha_{pre}}$ back into the original optimal guidance command $\alpha$ for practical implementation. It can be demonstrated that the sign of the original optimal guidance command $\alpha$ is always opposite to that of $\lambda_{v}$, i.e.,
\begin{align}
   4 \lambda_{v} \cdot \arctan \alpha   = -3\lambda_{u} - \sqrt{9\lambda_{u}^2 + 8\lambda_{v}^2} < 0. 
\end{align}
Hence, the original optimal guidance command $\alpha$ can be obtained by simply comparing the sign of the preprocessed guidance command $\alpha_{pre}$ and $\lambda_{v}$, i.e.,
\begin{align}
   \alpha = 
   \begin{cases}
      -\alpha_{pre}, \text{if}~\alpha_{pre}\cdot \lambda_{v} >0\\
      \alpha_{pre},~~\text{otherwise}
   \end{cases}
\end{align}
\subsection{Scheme for Generating Optimal Guidance Commands in Real Time}\label{SchemeGuidanceCommands}
\begin{figure}[!htp]
\centering
\includegraphics[width = 0.8\linewidth]{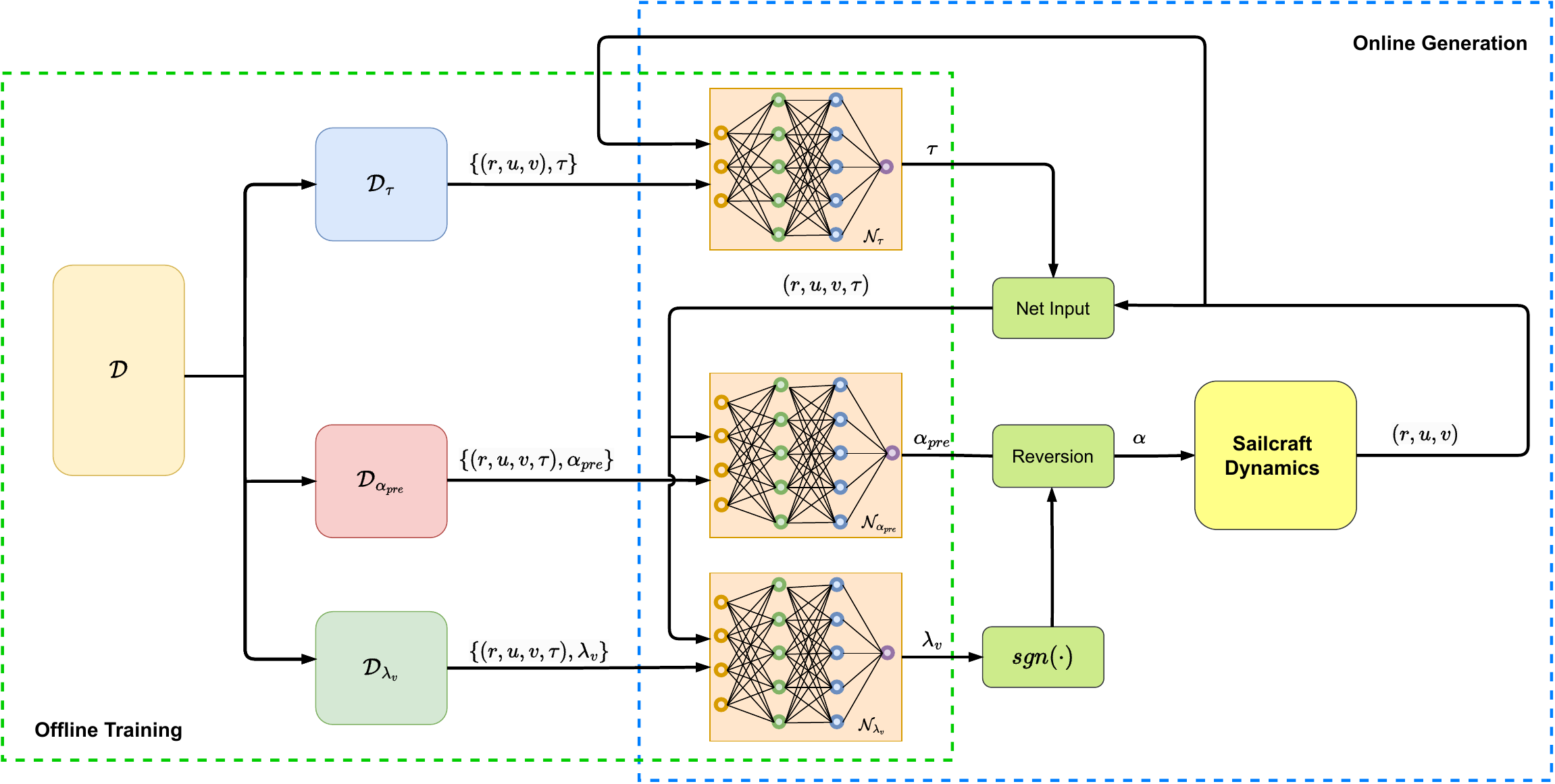}
\caption{Generation of optimal guidance commands in real time via NNs.}
\label{Fig:dnn_three}
\end{figure}

To facilitate real-time implementation, a system composed of three NNs is established, as depicted in Fig.~\ref{Fig:dnn_three}. To elaborate, $\mathcal{D}$ is divided into three samples, i.e., $\mathcal{D}_{\tau}:=\left\{(r,u,v),\tau \right\}$, $\mathcal{D}_{\alpha_{pre}}:=\left\{(r,u,v,\tau),\alpha_{pre} \right\}$, and $\mathcal{D}_{\lambda_{v}}:=\left\{(r,u,v,\tau), \lambda_{v} \right\}$; and the corresponding NN is denoted by $\mathcal{N}_{\tau}$, $\mathcal{N}_{\alpha_{pre}}$, and $\mathcal{N}_{\lambda_{v}}$, respectively. $\mathcal{N}_{\tau}$ is designed to forecast the optimal time of flight $\tau$ given a flight state $(r, u, v)$. 
Additionally, this network facilitates preliminary mission design by providing time of flight evaluations without necessitating an exact  solution.  
The output of $\mathcal{N}_{\tau}$ is subsequently used as part of the input for $\mathcal{N}_{\alpha_{pre}}$, which yields the preprocessed guidance command $\alpha_{pre}$. 
Concurrently, $\mathcal{N}_{\lambda_{v}}$ predicts the co-state $\lambda_{v}$, in which its sign $sgn(\lambda_{v})$ is used to revert the preprocessed guidance command $\alpha_{pre}$ back into the original guidance command $\alpha$. 
Once these three NNs are trained offline, they enable the generation of the optimal guidance command $\alpha$ given an initial condition $(r_0, u_0, v_0)$.
Consequently, the trained NNs offer closed-form solutions to the OCP. This endows the proposed method with robustness and generalization abilities, as shown in Subsections ~\ref{robus}  and \ref{gener}.
\subsection{NN Training}
Now, our focus shifts to the implementation of NN training algorithm. In addition to the three aforementioned NNs, i.e., $\mathcal{N}_{\tau}$, $\mathcal{N}_{\lambda_{v}}$, and $\mathcal{N}_{\alpha_{pre}}$, we also include the training of $\mathcal{N}_{\alpha}$ based on the sample $\mathcal{D}_{\alpha}=\left\{(r,u,v,\tau),\alpha \right\}$.
This is done to highlight the enhancement in approximation resulting from the preprocessing procedure. All the networks considered are feedforward NNs with multiple hidden layers. 
\it{Prior} \rm{to} training, the dataset samples \textcolor{black}{are shuffled} to establish a split of $70\%$ for training, $15\%$ for validation, and $15\%$ for testing sets. \textcolor{black}{All input and output data are normalized by the min-max scaling and scaled to the range [0, 1].} Selecting an appropriate NN structure, particularly the number of hidden layers and neurons per layer, is a non-trivial task. A structure that is too simplistic tends to lead to underfitting. Conversely, overfitting often arises when the structure is overly complex. Making a balance between time consumption and the overfitting issue, we adopt a structure with three hidden layers, each containing 20 neurons. Subsequently, the sigmoid function serves as the activation function. For the output layer, a linear function is utilized.
The crux of the training lies in minimizing the loss function, quantified as the mean squared error (MSE) between the predicted values from the trained NNs and the actual values within the dataset samples. We employ the 'Levenberg-Marquardt' in Ref. \cite{hagan1994training} for training NNs, \textcolor{black}{and the training is terminated after 1,000 epochs or when the loss function drops below $1\times 10^{-8}$}. To prevent overfitting, we employ an early stopping criteria based on the concept of patience. Other hyperparameters utilized in training are set to their default values.

\begin{figure}[!htp]
\centering
\begin{subfigure}[t]{7cm}
\centering
\includegraphics[width = 8cm]{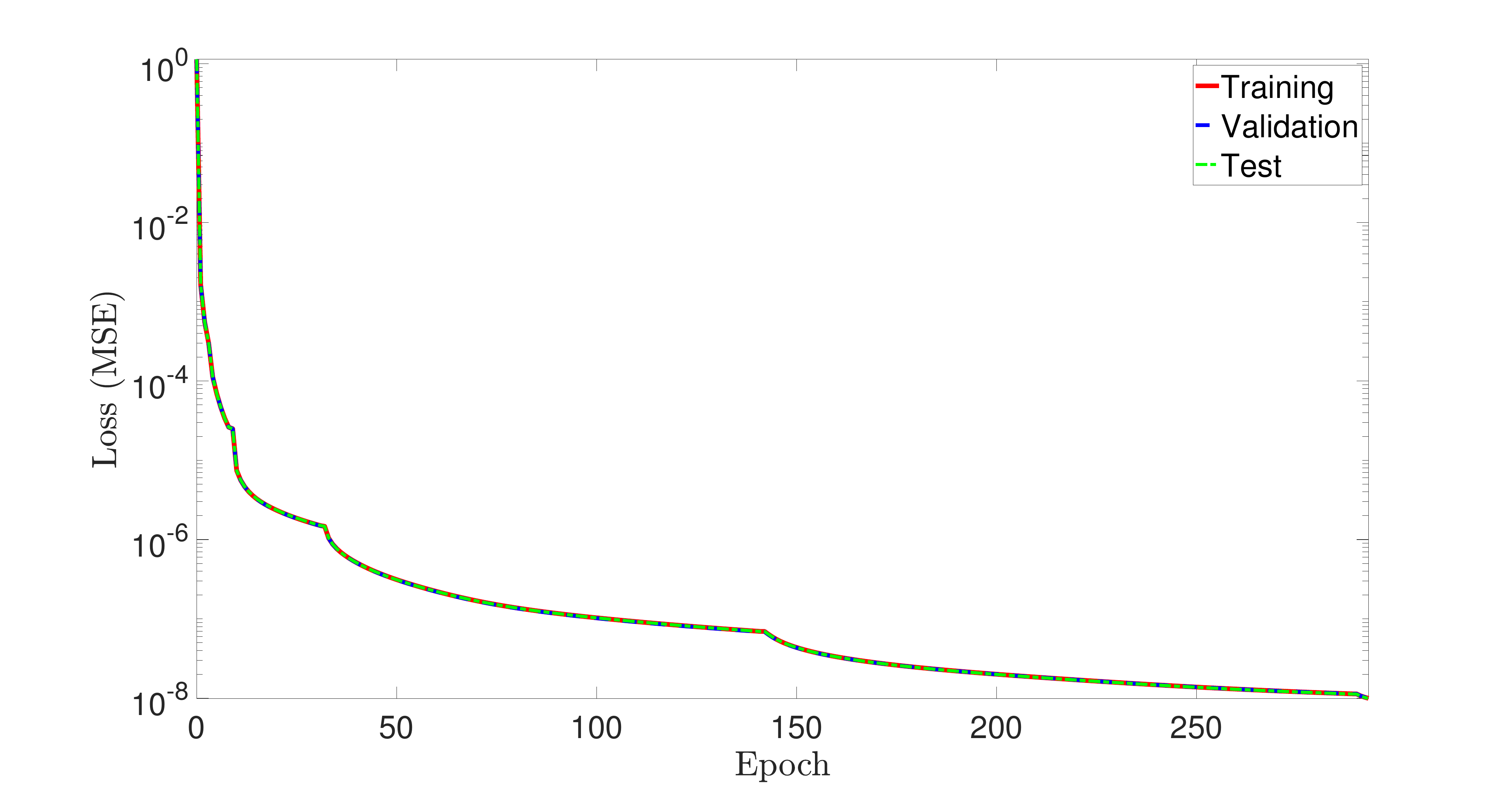}
\caption{$\mathcal{N}_{\tau}$}
\label{Fig:loss_tf}
\end{subfigure}
~~~~~
\begin{subfigure}[t]{7cm}
\centering
\includegraphics[width = 8cm]{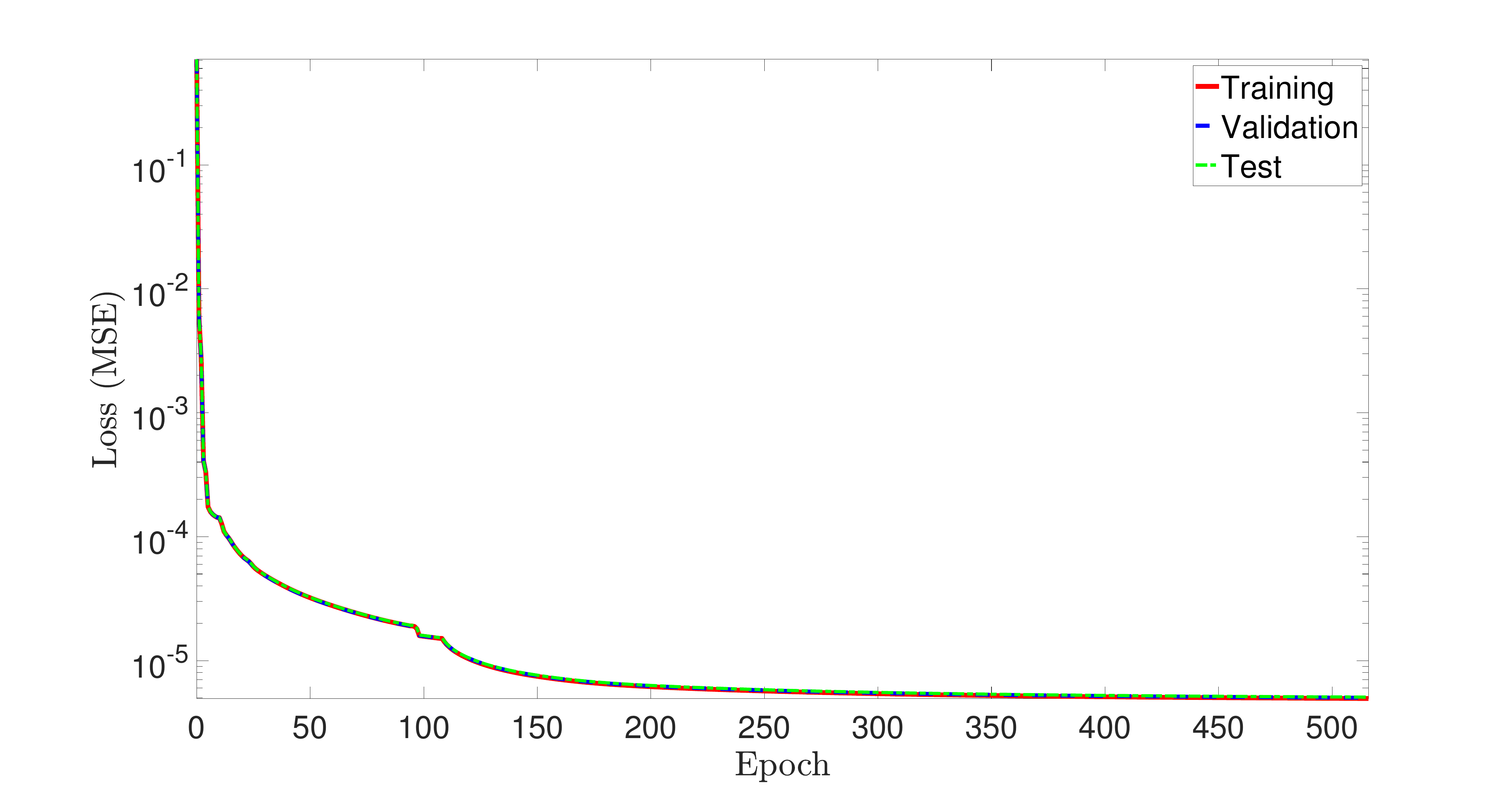}
\caption{$\mathcal{N}_{\lambda_{v}}$}
\label{Fig:loss_pvtheta}
\end{subfigure}\\
\begin{subfigure}[t]{7cm}
\centering
\includegraphics[width = 8cm]{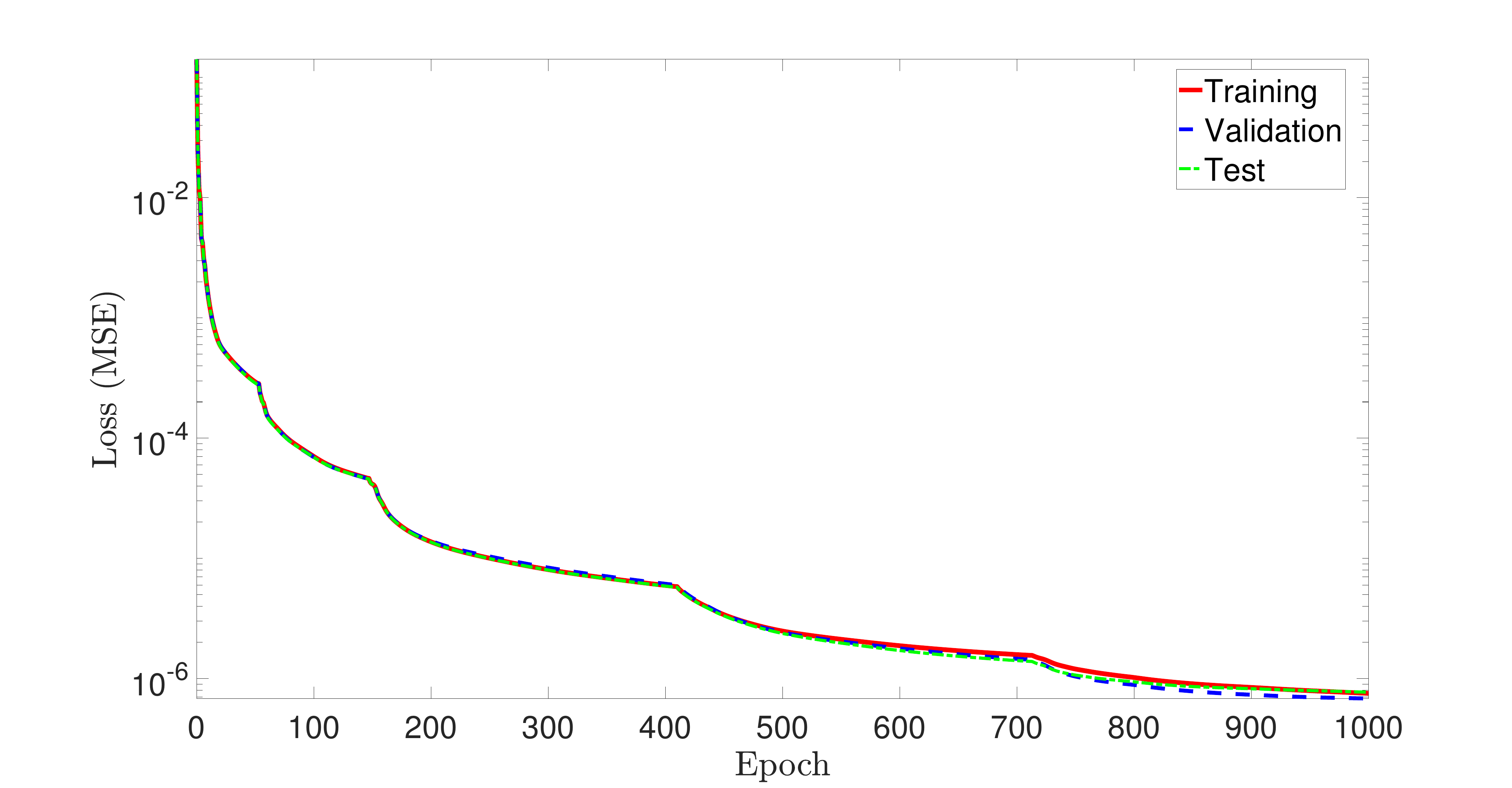}
\caption{$\mathcal{N}_{\alpha_{pre}}$}
\label{Fig:loss_reg_u}
\end{subfigure}
~~~~~
\begin{subfigure}[t]{7cm}
\centering
\includegraphics[width = 8cm]{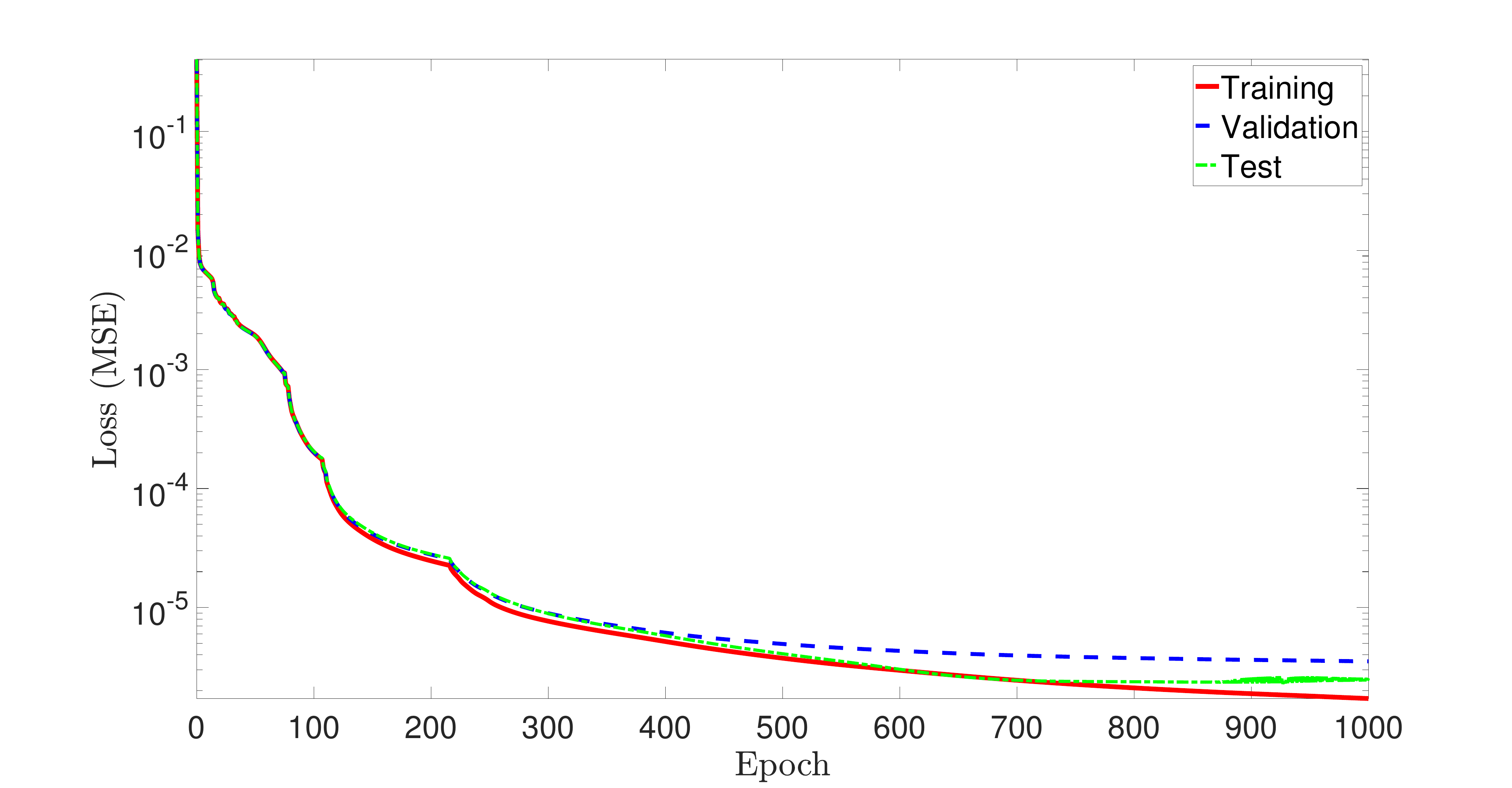}
\caption{$\mathcal{N}_{\alpha}$}
\label{Fig:loss_alpha}
\end{subfigure}
\caption{Training histories of four NNs.}
\label{Fig:Training}
\end{figure}

Fig.~\ref{Fig:Training} 
illustrates the training progression of four NNs. As seen in Fig.~\ref{Fig:loss_tf} , it is evident that the MSEs for the training, validation, and test sets decrease to $1\times 10^{-8}$ in less than $300$ epochs. Table~\ref{Table:control_effort_25_50} displays the MSEs of the four NNs upon completion of the training process. Notably, the validation and test errors of $\mathcal{N}_{\alpha_{pre}}$ are smaller than those of $\mathcal{N}_\alpha$, indicating that our preprocessing procedure enhances the approximation accuracy. This is further confirmed through simulations detailed in the next section.
\begin{table}[!htp]
\centering
\caption{MSEs of four NNs}
\begin{tabular}{ccccc}
\hline
   & $\mathcal N_{\tau}$      & $\mathcal{N}_{\lambda_{v}}$    & $\mathcal{N}_{\alpha_{pre}}$  & $\mathcal{N}_{\alpha}$\\
\hline
Training &$1\times 10^{-8}$ &$ 4.94\times 10^{-6}$ &$7.53 \times 10^{-7}$ &$ 1.72\times 10^{-6}$ \\
Validation &$1\times 10^{-8}$ &$ 5.02\times 10^{-6}$  &$6.82\times 10^{-7}$ &$ 3.53\times 10^{-6}$ \\ 
Test &$1\times 10^{-8}$ &$5.07\times 10^{-6}$ &$7.70\times 10^{-7}$ &$2.55\times 10^{-6}$ \\ 
\hline
\label{Table:control_effort_25_50}
\end{tabular}
\end{table}
\section{Numerical Simulations}
\label{NumericalSimulations}
In this section, we present some numerical simulations to demonstrate the efficacy and performance of the proposed methodology. We commence by showing its optimality performance. Subsequently, robustness against perturbations and generalization ability of the proposed method are investigated. In addition, we evaluate the real-time execution of the guidance command generation. It is worth noting that the proposed method has been implemented on a desktop equipped with an Intel Core i9-10980XE CPU @3.00 GHz and 128 GB of RAM.
\subsection{Optimality Performance}
\label{firstSim}
This subsection is devoted to examining whether or not the proposed approach can produce optimal solutions compared with existing methods. 
Table~\ref{Table:Initialconditions} outlines the initial conditions for two orbital transfers. For comparison, two strategies are employed to steer the solar sailcraft. The first strategy employs three NNs, i.e., $\mathcal{N}_{\tau}$, $\mathcal{N}_{\alpha_{pre}}$, and $\mathcal{N}_{\lambda_{v}}$, while the second one uses two NNs, $\mathcal{N}_{\tau}$ and $\mathcal{N}_{\alpha}$. The indirect method, resolving the shooting function in Eq.~(\ref{EQ:TPBVP_law}) will be used for comparison.
\begin{table}[!htp]
\centering
\caption{Initial conditions for two orbital transfer cases}
\begin{tabular}{ccccc}
\hline
   & $r_0$      & $\theta_0$    & $u_0$  & $v_0$\\
\hline
Case 1 &$1$ &$0$ &$0$ &$1$ \\
Case 2 &$1.05$ &$0$  &$0.15$ &$1.03$ \\ 
\hline
\label{Table:Initialconditions}
\end{tabular}
\end{table}

\begin{figure}[!htp]
\centering
\begin{subfigure}[t]{7cm}
\centering
\includegraphics[width = 8cm]{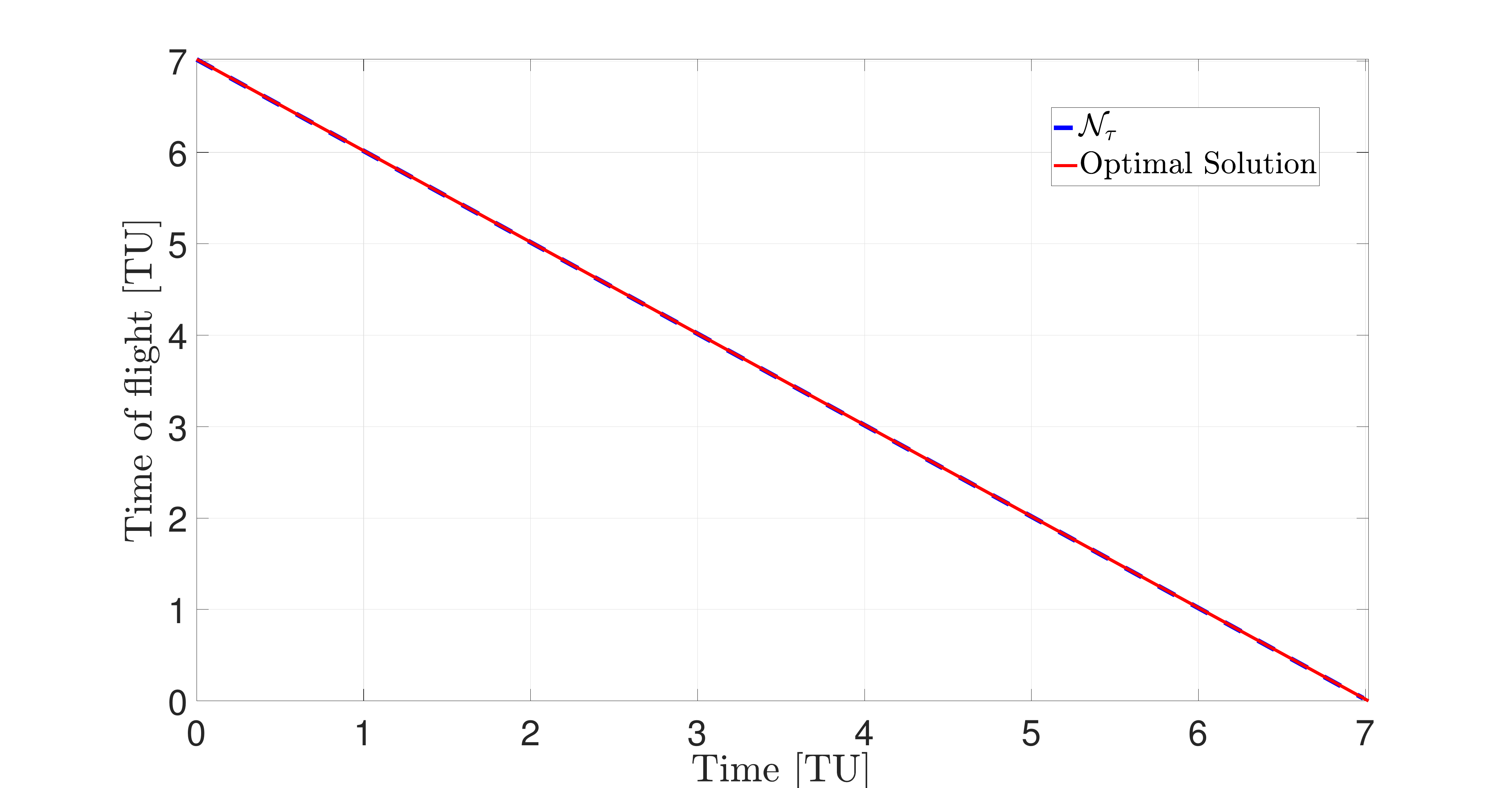}
\caption{Time of flight profile}
\label{Fig:Time_of_flight}
\end{subfigure}
~~~~~
\begin{subfigure}[t]{7cm}
\centering
\includegraphics[width = 8cm]{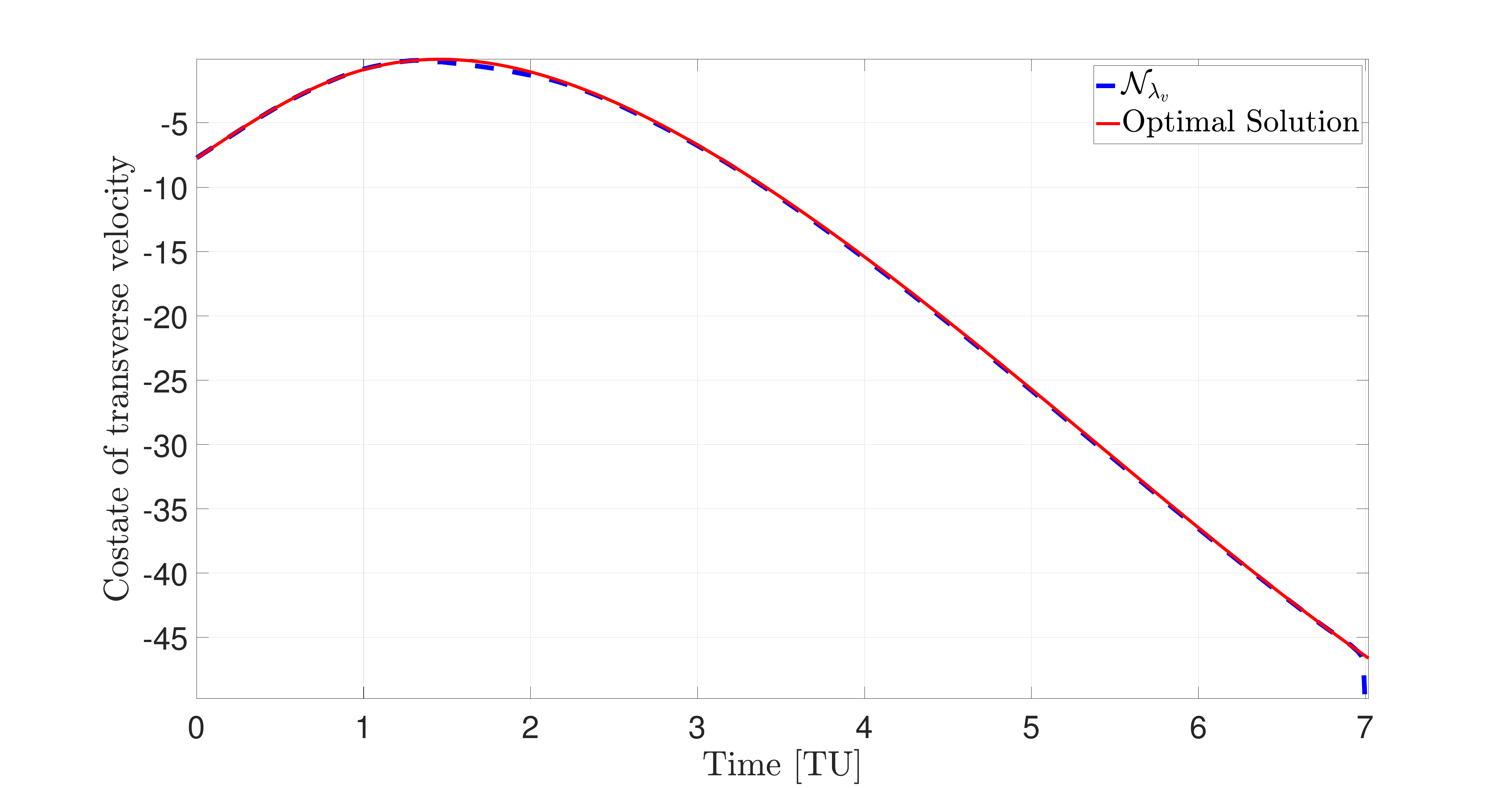}
\caption{Profile for co-state of the transverse speed}
\label{Fig:Costate_oftangential}
\end{subfigure}\\
\begin{subfigure}[t]{7cm}
\centering
\includegraphics[width = 8cm]{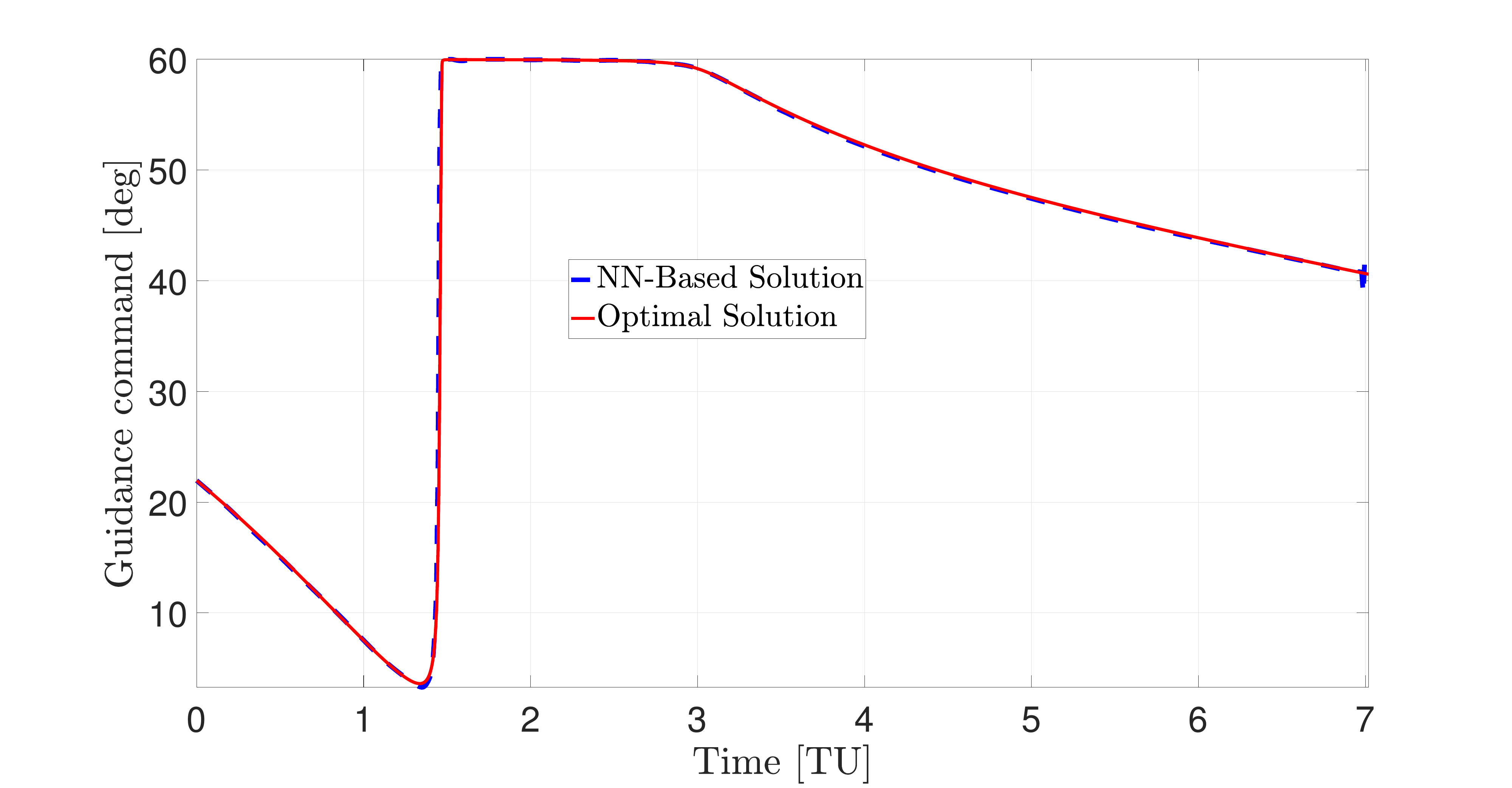}
\caption{Guidance command profile derived from $\mathcal{N}_{\alpha_{pre}}$, $\mathcal{N}_{\lambda_{v}}$ and $\mathcal{N}_{\tau}$.}
\label{Fig:cooperativesteering}
\end{subfigure}
~~~~~
\begin{subfigure}[t]{7cm}
\centering
\includegraphics[width = 8cm]{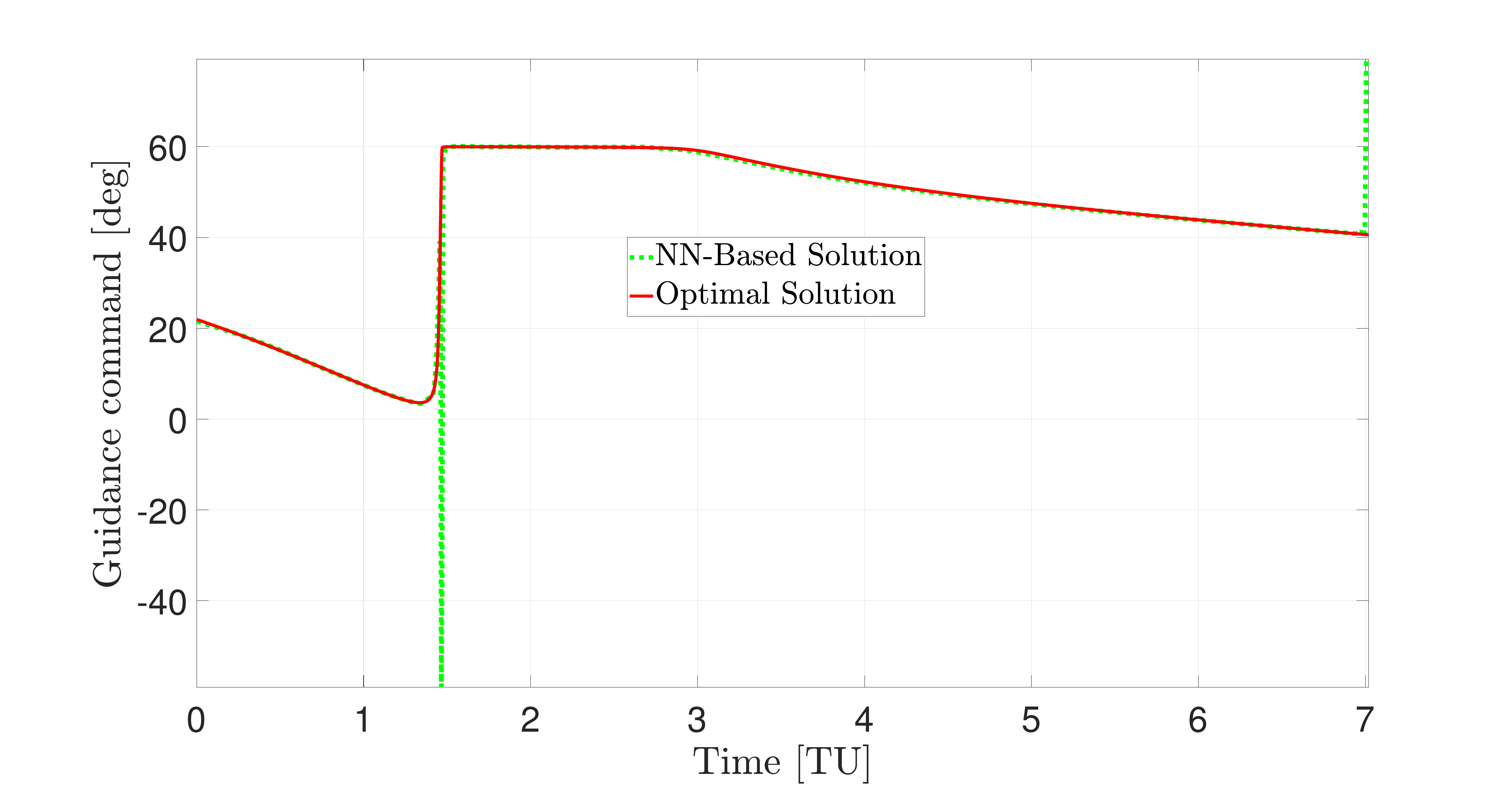}
\caption{Guidance command profile derived  from $\mathcal{N}_{\alpha}$ and $\mathcal{N}_{\tau}$.}
\label{Fig:cooperativeVelocity}
\end{subfigure}
\caption{Comparison of results from the trained NNs and the indirect method for Case 1.}
\label{Fig:cooperative_control_profile}
\end{figure}

For Case 1, the solar sailcraft embarks on a journey from Earth's orbit around the Sun, and the mission is to steer the solar sailcraft into Mars' orbit around the Sun. 
Fig.~\ref{Fig:cooperative_control_profile} represents the outcomes yielded by the trained NNs and the indirect method. The optimal time of flight obtained via the indirect method is $7.0184$ TU, while the prediction from $\mathcal{N}_{\tau}$ is $7.0192$ TU for the given initial condition. 
Remarkably, $\mathcal{N}_{\tau}$ exhibits precise approximation of the optimal time of flight throughout the entire transfer, as shown by Fig.~\ref{Fig:Time_of_flight}. On the other hand, Fig.~\ref{Fig:Costate_oftangential} shows less accurate approximation of the optimal co-state $\lambda_{v}$, especially at the end of the transfer. 
Fortunately, the sign of the predicted $\lambda_{v}$, rather than its exact value, is sufficient for reverting the preprocessed guidance command back to the original guidance command. This observation is depicted in Fig.~\ref{Fig:cooperativesteering}, where the guidance commands derived from the first strategy \textcolor{black}{closely align with the indirect method, except at the end of the transfer. Since the chattering guidance command derived from the proposed method at the end of the transfer is not practically implementable, the numerical simulation is terminated once the predicted time of flight $\tau$ from $\mathcal{N}_{\tau}$ drops below 0.005 TU (equivalent to 0.2907 days) hereafter.}
In contrast, the guidance command derived from the second strategy exhibits lower accuracy, especially in the presence of rapid changes in the optimal guidance command. The degeneration of approximation performance, caused by the second strategy, is demonstrated in Fig.~\ref{Fig:cooperativeVelocity}. It is worth emphasizing that this problem of performance degeneration may not be well fixed by simply adjusting the hyperparameters for training NNs, as shown by the numerical simulations in Ref. \cite{li2019neural}. In that work, a random search for the hyperparameters was employed, and the trained NNs, even with the best approximation performance, were not able to capture the rapid change in the guidance command. In contrast, thanks to the preprocessing procedure in this paper, the well-trained NNs are capable of precisely generating the optimal guidance command.

To assess the performance of fulfilling the terminal condition, we define the flight error $\Delta \Phi$ as the difference between the NN-driven terminal flight state $(r^{\mathcal{N}}(t_f),u^{\mathcal{N}}(t_f),v^{\mathcal{N}}(t_f))$ and desired terminal flight state $(1.524,0,\frac{1}{\sqrt{1.524}})$. This error is computed as follows
\begin{align*}
\Delta \Phi = ( \lvert r^{\mathcal{N}}{(t_f)}-1.524 \rvert,\lvert u^{\mathcal{N}}(t_f)\rvert,\lvert v^{\mathcal{N}}(t_f)-\frac{1}{\sqrt{1.524}}\rvert).
\end{align*}

\begin{figure}[!htp]
\centering
\includegraphics[width = 0.7\linewidth]{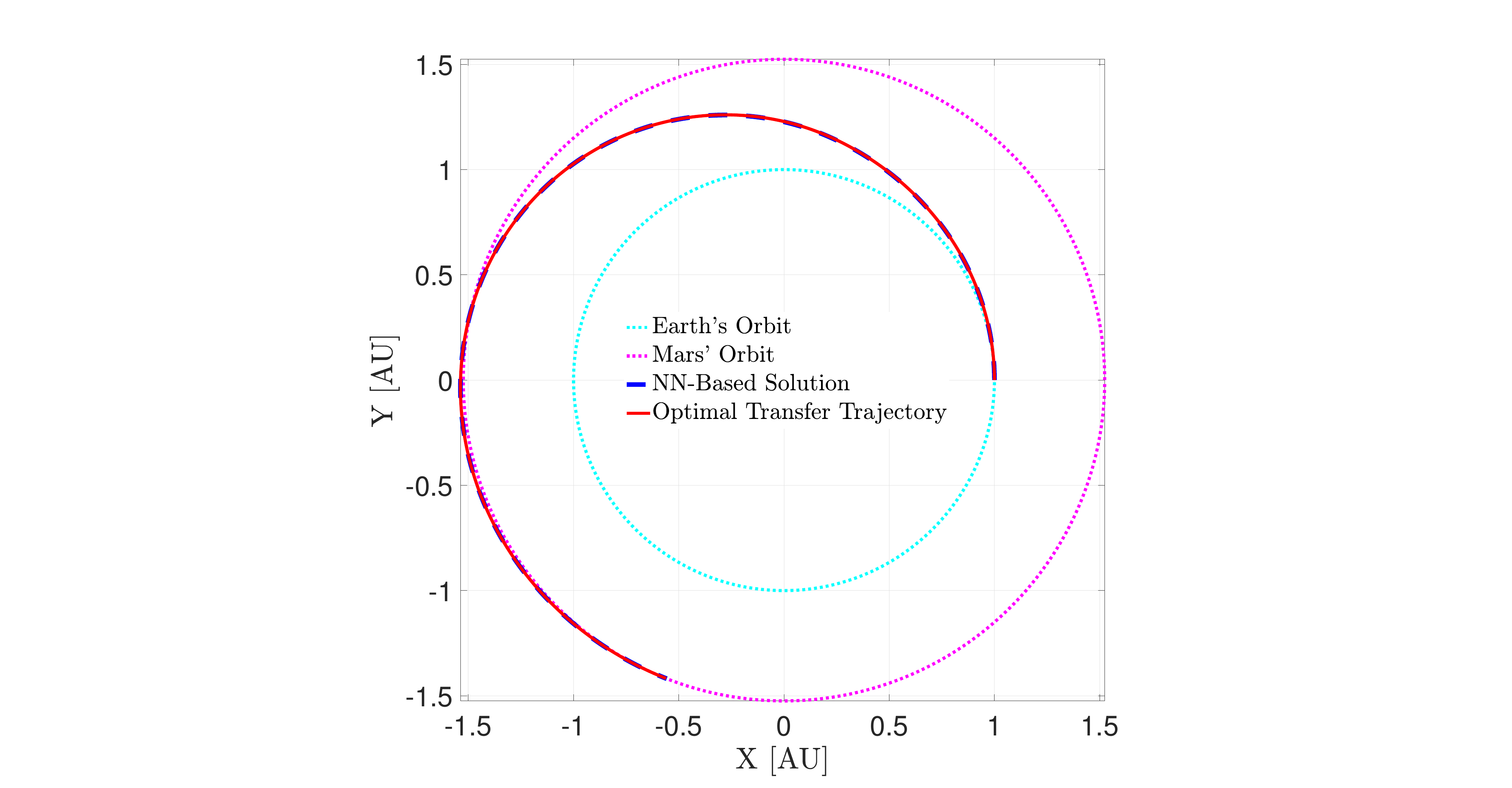}
\caption{Transfer trajectories derived from the first strategy and the indirect method for Case 1.}
\label{Fig:tra_nominalcase1}
\end{figure}

Then, the flight error caused by the first strategy is $\Delta \Phi = (5.3703 \times 10^{-7}, 3.0121 \times 10^{-5}, 3.6943 \times 10^{-5} )$, which is lower than that of the second strategy with $\Delta \Phi =( 2.2220 \times 10^{-5}, 4.1852 \times 10^{-4}, 2.5968 \times 10^{-4})$. Although both strategies can guide the solar sailcraft into Mars' orbit with negligible flight errors, 
the guidance command derived from the first strategy is more feasible and accurate than the second strategy, as shown in Figs.~\ref{Fig:cooperativesteering} and \ref{Fig:cooperativeVelocity}. 
In addition, the transfer trajectories derived from the first strategy and the indirect method are illustrated in Fig.~\ref{Fig:tra_nominalcase1}. 

\begin{figure}[!htp]
\centering
\includegraphics[width = 0.7\linewidth]{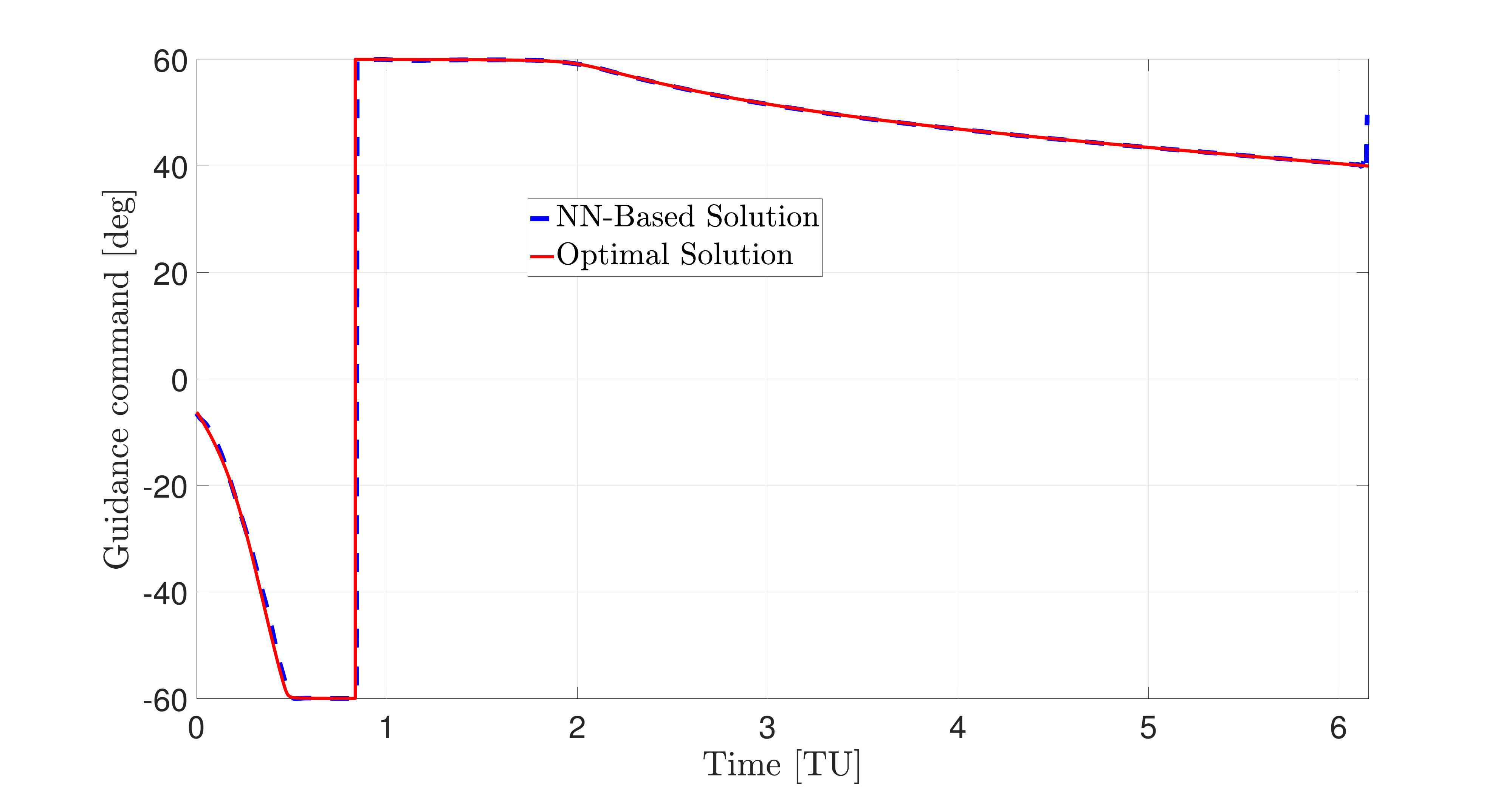}
\caption{Guidance command profiles derived from the NN-based approach
and the indirect method for Case 2.}
\label{Fig:control_nominalcase2}
\end{figure}

\begin{figure}[!htp]
\centering
\includegraphics[width = 0.7\linewidth]{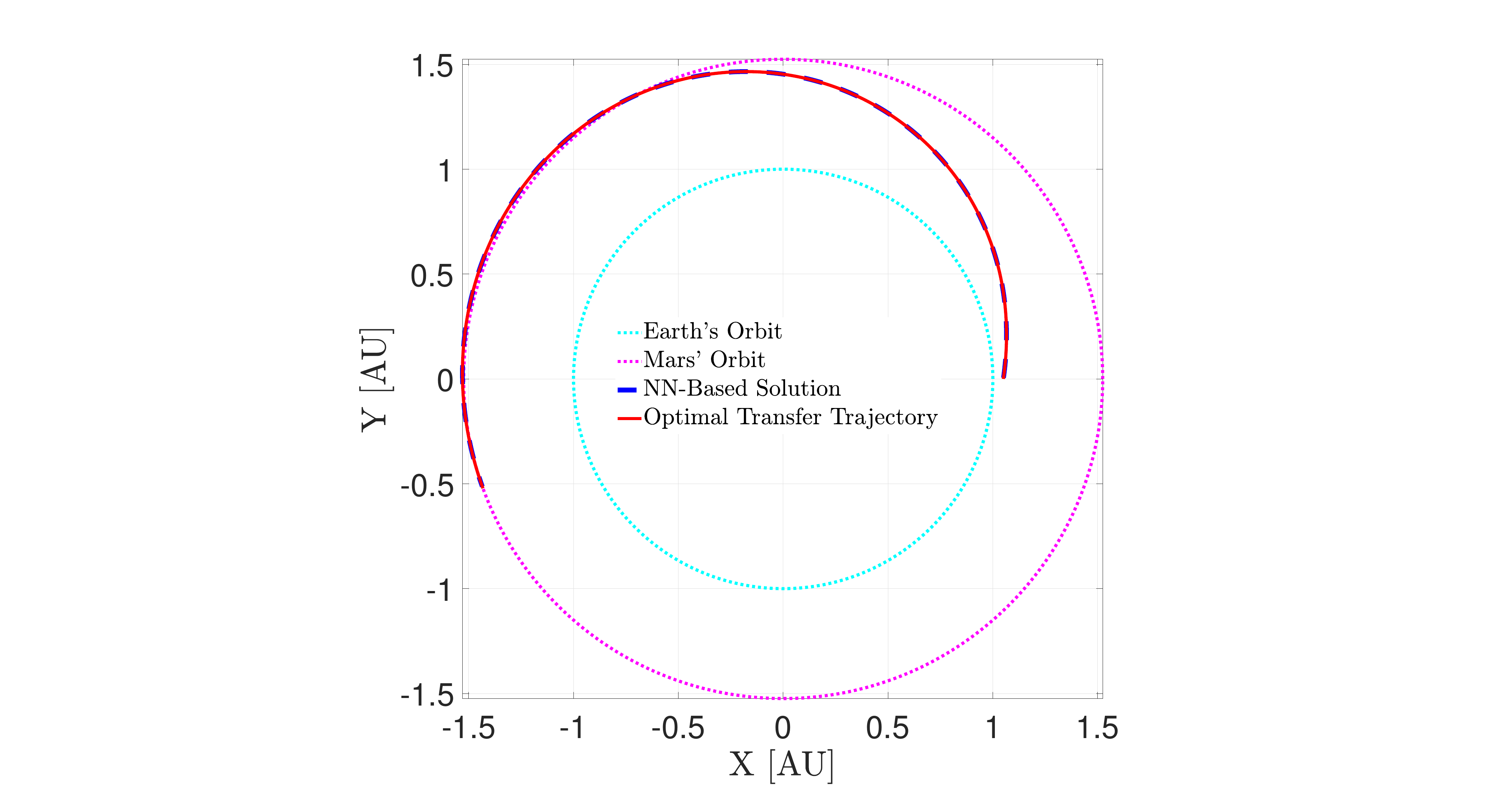}
\caption{Transfer trajectories derived from the NN-based approach and the indirect method for Case 2.}
\label{Fig:transfer_nominal}
\end{figure}

Given the superiority of the first strategy over the second, our focus will exclusively be on the first strategy, which will be called NN-based approach hereafter. 
Shifting to Case 2, the output of $\mathcal{N}_{\tau}$ for the initial condition stands at 6.1643 TU, while the optimal solution attained via the indirect method is 6.1552 TU, resulting in a marginal functional penalty of merely 0.5290 days. Furthermore, 
Fig.~\ref{Fig:control_nominalcase2} illustrates that the NN-based approach precisely captures the discontinuous jump and effectively approximates the optimal guidance command, except at the end of the transfer. The transfer trajectories derived from the NN-based approach and the indirect method are depicted in Fig.~\ref{Fig:transfer_nominal}.
As a consequence, the flight error caused by the NN-based approach is $\Delta \Phi = (5.2011 \times 10^{-6}, 3.3526 \times 10^{-4}, 1.5393 \times 10^{-4})$, implying that the NN-based approach capably steers the solar sailcraft into Mars' orbit with \textcolor{black}{an acceptable flight error.} 
\subsection{Robustness Analysis}\label{robus}
In reality, the normalized characteristic acceleration $\beta$ is influenced by diverse factors, such as solar sailing efficiency \textcolor{black}{and} solar radiation pressure \cite{spencer2019solar}. Consequently, the solar sailcraft must continually adjust its flight trajectory based on its current flight state to rectify any flight errors. Unfortunately, this process is quite challenging for onboard systems with limited computational resources. Additionally, solving OCPs featuring parameters under perturbations tends to be problematic for conventional direct and indirect methods that rely on precise system models. For this reason, we evaluate the proposed method's robustness against perturbations concerning $\beta$.

The initial condition for the solar sailcraft is set as
\begin{align*}
r_0 = 1.1, \theta_0 = \frac{\pi}{2}, u_0 = 0.18, v_0 = 0.93,
\end{align*}
and the mission is to fly into Mars' orbit. The normalized characteristic acceleration $\beta$ is under perturbations following a standard uniform distribution within the range of $(-15\%, 15\%)$. 

\begin{figure}[!htp]
\centering
\includegraphics[width = 0.7\linewidth]{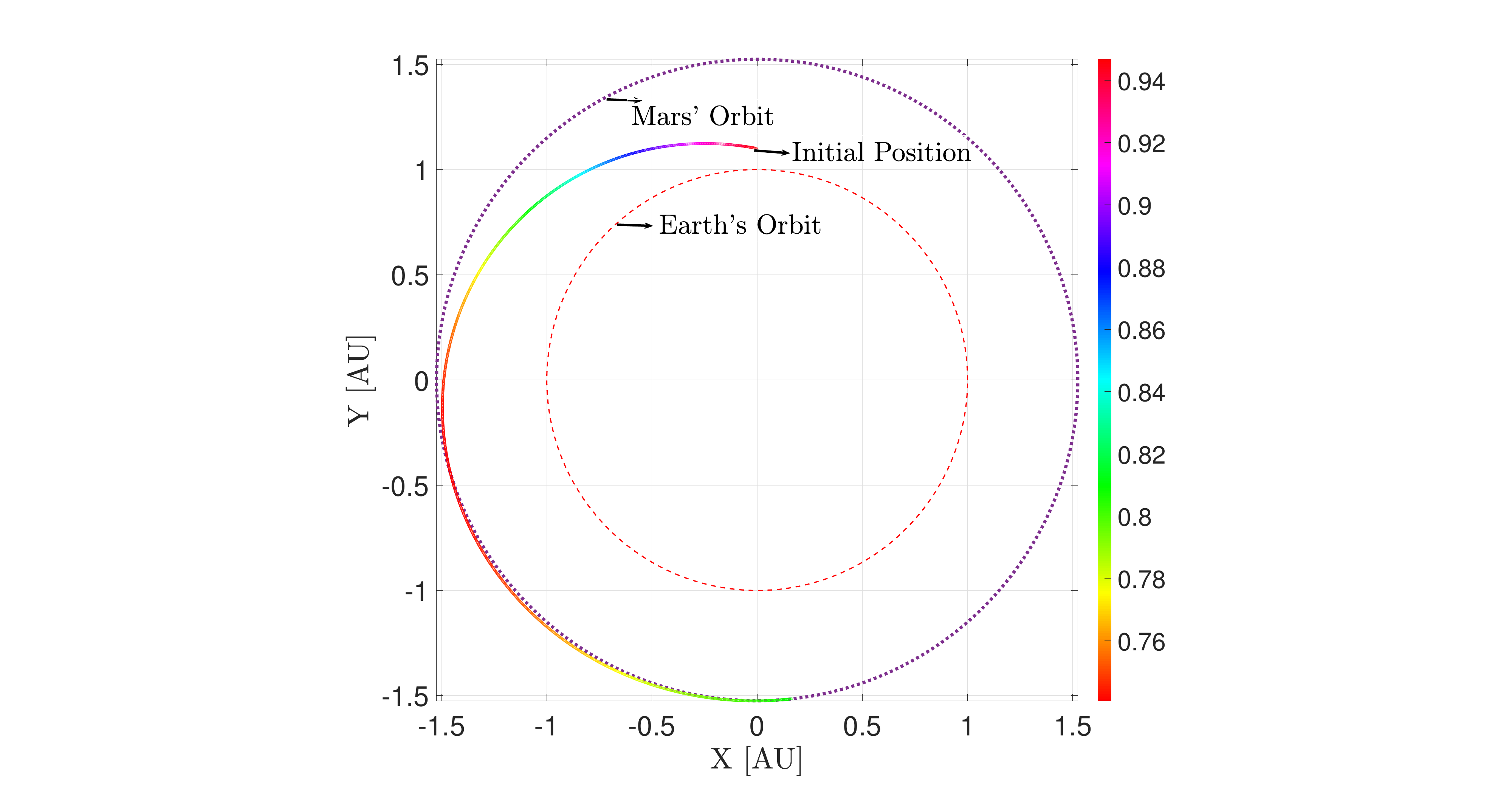}
\caption{NN-based transfer trajectory with $\beta$ under perturbations.}
\label{Fig:transfer_noise}
\end{figure}

\begin{figure}[!htp]
\centering
\includegraphics[width = 0.7\linewidth]{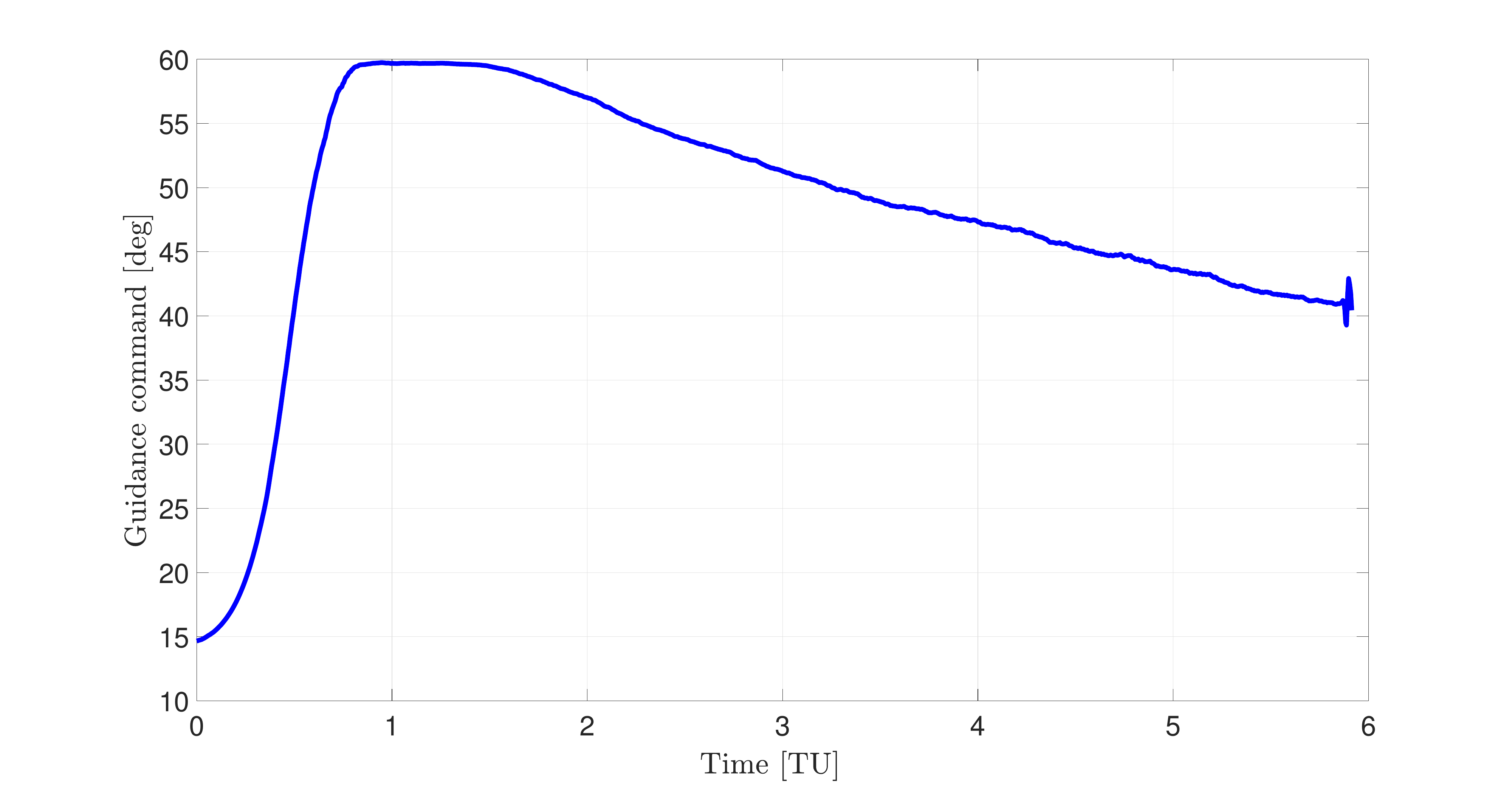}
\caption{NN-based guidance command profile with $\beta$ under perturbations.}
\label{Fig:control_noise}
\end{figure}

Figs.\ref{Fig:transfer_noise} and \ref{Fig:control_noise} depicts the transfer trajectory and the corresponding guidance command derived from the NN-based approach, respectively.  
It can be seen that the colormap-identified speed exhibits a gradual reduction for most of the transfer, followed by a gradual increase
towards the end of the transfer, as depicted in Fig.\ref{Fig:transfer_noise}. Regarding Fig.\ref{Fig:control_noise}, the guidance command is generally smooth during the initial phase of the transfer.
Subsequently, to uphold precise flight, the NN-based approach acts like an error corrector by dynamically adapting the guidance command in the latter phase of the transfer. A remarkable observation is that even when the solar sailcraft's dynamics is under perturbations, the attitude constraints remain unviolated throughout the entire transfer. Furthermore, it takes 5.9152 TU for the solar sailcraft to accomplish its journey with a minor flight error of $\Delta \Phi = (2.0160 \times 10^{-5}, 9.8708 \times 10^{-5}, 8.4738 \times 10^{-5})$.
\subsection{Generalization Analysis}\label{gener}
\begin{figure}[!htp]
\centering
\includegraphics[width = 0.7\linewidth]{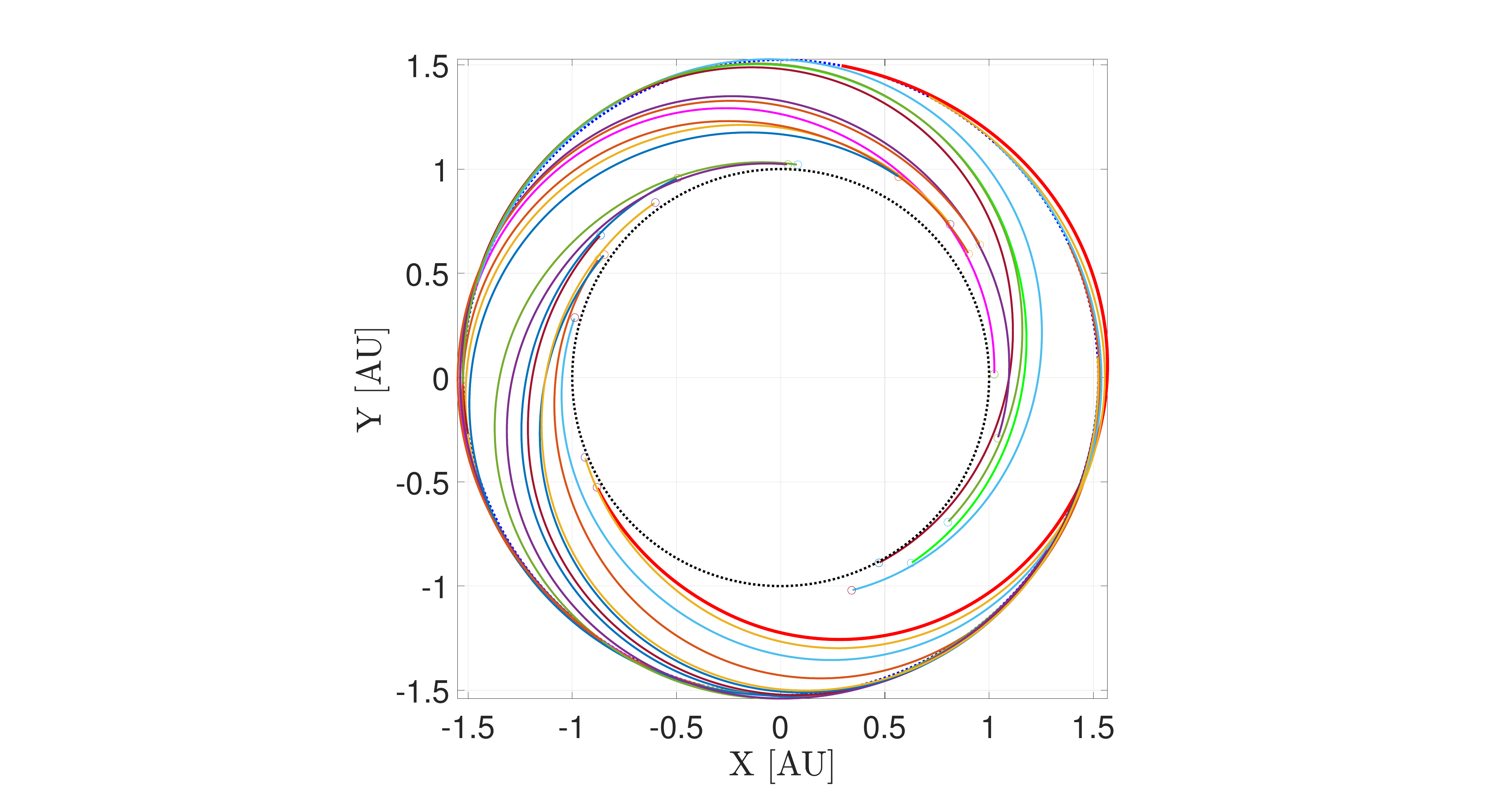}
\caption{NN-based transfer trajectories with initial conditions randomly chosen in $\mathcal{A}$.}
\label{Fig:generalisation}
\end{figure}

The assessment of generalization ability pertains to evaluating the predictive performance for inputs that lie outside the training dataset. Define a state space $\mathcal{A}$, that is not inside  $\mathcal{D}$, as below
\begin{align*}
r_0 \in [1,1.15], \theta_0 \in [0, 2\pi], u_0 \in [0,0.1], v_0 \in [0.8,1.2].
\end{align*}
Within this space, the solar sailcraft's initial conditions are randomly chosen. Subsequently, the proposed approach is employed to guide the solar sailcraft from the selected initial condition to Mars' orbit. A total of 30 tests are conducted, resulting in 20 successful transfers, as depicted in Fig.~\ref{Fig:generalisation}, in which the small circles denote the solar sailcraft's initial positions. Among these successful cases, the average flight error is $\Delta \Phi = ( 4.5091 \times 10^{-5}, 1.4726 \times 10^{-4},1.0147 \times 10^{-4})$. This indicates that the proposed method is able to steer the spacecraft to Mars’ obit accurately \textcolor{black}{even with some initial states that are not in the training set.} 
\subsection{Real-Time Performance}
Since the output of a feedforward NN is essentially a composition of linear mappings of the input vector, a trained NN can produce the output within a constant time. Recall that our method requires three trained NNs with simple structures. 10,000 trials of the proposed method across various flight states are run in a C-based computational environment, and the mean execution time for generating a guidance command is 0.0177 ms. This translates to approximately 0.5310 ms on a typical flight processor operating at 100 MHz \cite{gankidi2017fpga}. Conversely, indirect methods typically necessitate time-consuming processes of integration and iteration. We utilize \it{fsolve} \rm{function} to solve the relevant TPBVP and observe a convergence time of approximately 1.22 seconds, even when an appropriate initial guess is provided. Consequently, the computational burden associated with optimization-based indirect methods can be substantially alleviated through adopting the proposed method.
\section{Conclusions} \label{Colus}
Real-time optimal control for attitude-constrained solar sailcrafts in coplanar circular-to-circular interplanetary transfers was studied in the paper. The necessary conditions governing optimal trajectories were deduced through the PMP. Different from the conventional methods that rely on optimization-based solvers to \textcolor{black}{collect} a training dataset, we formulated a parameterized system. This enabled us to generate an optimal trajectory by simply solving an initial value problem, and the optimal state-guidance command pairs could be readily obtained. Regarding the challenge posed by potential discontinuities in the optimal guidance \textcolor{black}{command}, which could undermine the NN's \textcolor{black}{approximation} performance, we developed a technique to preprocess the guidance command. With the saturation function introduced in the optimal guidance law, this preprocessing technique smoothened the guidance command. \textcolor{black}{In this way, the approximation performance of the preprocessed guidance command was enhanced, and another NN was used to predict one co-state, whose sign was adopted to revert the preprocessed guidance command into the potentially discontinuous original guidance command.
As a result, the well-trained NNs were capcable of 
not only generating the optimal guidance command precisely in real time but also predicting the optimal time of flight given a flight state.} 
\section*{Acknowledgement}
This research was supported by the National Natural Science Foundation of China under Grant No. 62088101. 

\appendix
\setcounter{equation}{0}
\renewcommand\theequation{A\arabic{equation}} 
\renewcommand\thefigure{\Alph{section}\arabic{figure}} 
\section{Proof \textcolor{black}{of} Lemmas in Section \ref{Problem}}\label{Appendix:A}

Proof \textcolor{black}{of} Lemma \ref{LE:optimal_control_law}. 
In view of $\dot{\lambda}_{\theta}(t) = 0$ in Eqs.~(\ref{EQ:costate_function}) and (\ref{EQ:theta_law}), it is easy to see that $\lambda_{\theta}(t)$ remains zero for $t \in [0,t_f]$ along an optimal trajectory. Then,  Eq.~(\ref{EQ:Ham_function}) reduces to
\begin{align}
\mathscr H = 1 + \lambda_r u  + \lambda_{u}(\frac{\beta \cos^3\alpha}{r^2} + \frac{v^2}{r} - \frac{1}{r^2}) + \lambda_{v}
(\frac{\beta \sin\alpha \cos^2\alpha}{r^2} - \frac{u v}{r}).
\label{A1_H}
\end{align}
To obtain stationary values of $\mathscr H$, we have
\begin{align}
\frac{\partial \mathscr H}{\partial \alpha} = 
\frac{\beta}{r^2}\cos ^3 \alpha [-3\lambda_{u} \tan \alpha + \lambda_{v}(1-2 \tan ^2 \alpha)] = 0.
\label{A1_stationary}
\end{align}
Assume that $\lambda_{u}$ and $\lambda_{v}$ are not both zero \cite{wood1982comment}. Because $\beta > 0$, and $\cos \alpha \neq 0$ for $\alpha \in [-\phi_{max},\phi_{max}]$, two different roots, denoted by $\alpha_1$ and $\alpha_2$,  can be obtained by solving Eq.~(\ref{A1_stationary}) as
\begin{align}
\begin{cases}
\alpha_{1} = \arctan ~\frac{-3\lambda_{u} + \sqrt{9 \lambda^2_{u} + 8 \lambda^2_{v}}}{4 \lambda_{v}},\\
\alpha_{2} = \arctan ~ \frac{-3\lambda_{u} - \sqrt{9 \lambda^2_{u} + 8 \lambda^2_{v}}}{4 \lambda_{v}}.
\label{A1_alpha}
\end{cases}
\end{align}
Now we proceed to check for the Hessian’s positive definiteness. According to Eq.~(\ref{A1_stationary}), taking the second partial derivative of $\mathscr H$ w.r.t. $\alpha$ yields
\begin{align}
\frac{\partial^2 \mathscr H}{\partial \alpha^2} = 
\frac{\beta}{r^2}\{-3 \cos^2\alpha \sin \alpha[-3\lambda_{u} \tan \alpha + \lambda_{v}(1-2 \tan^2 \alpha)] + \cos^3 \alpha[-3 \lambda_{u} \sec^2 \alpha - 4 \lambda_{v} \tan \alpha \sec^2\alpha]\}.
\label{A1_Hess}
\end{align}
Substituting Eq.~(\ref{A1_alpha}) into Eq.~(\ref{A1_Hess}) leads to
\begin{align}
\begin{cases}
\frac{\partial^2 \mathscr H}{\partial \alpha_1^2} = \frac{\beta}{r^2} \cos \alpha (-\sqrt{9 \lambda^2_{u} + 8 \lambda^2_{v}}),\\
\frac{\partial^2 \mathscr H}{\partial \alpha_2^2} = \frac{\beta}{r^2} \cos \alpha (+\sqrt{9 \lambda^2_{u} + 8 \lambda^2_{v}}).
\label{A1_Hess_alpha12}
\end{cases}
\end{align}
Clearly, $\frac{\partial^2 \mathscr H}{\partial \alpha^2} > 0$ \textcolor{black}{always} holds \textcolor{black}{if and only if} $\alpha = \alpha_2$. Therefore, the local minimum solution $\bar{\alpha}$ is given by
\begin{align}
\bar{\alpha} = \alpha_2 = \arctan~\frac{-3\lambda_{u} - \sqrt{9 \lambda^2_{u} + 8 \lambda^2_{v}}}{4 \lambda_{v}}.
\label{Eq:optimal_alpha}
\end{align}

Notice that Eq.~(\ref{Eq:optimal_alpha}) will become singular if $\lambda_{v}(t) = 0$. Then \textcolor{black}{$\lambda_{v}(t) = 0$ can hold at some isolated points for $t \in [0,t_f]$ or $\lambda_{v}(t) \equiv 0$ holds in a time interval $[t_a,t_b] \in [0,t_f]$}. We now analyze the first case. 

If $\lambda_{u}(t) < 0$ for $ t\in [t_a,t_b]$, then
\begin{align}
\begin{gathered}
\lim_{\lambda_{u} < 0, \lambda_{v} \to 0} \frac{-3\lambda_{u} - \sqrt{9 \lambda^2_{u} + 8 \lambda^2_{v}}}{4 \lambda_{v}} = \lim_{\lambda_{u} < 0, \lambda_{v} \to 0} \frac{-3\lambda_{u} +3\lambda_{u} \sqrt{1+ \frac{8 \lambda^2_{v}}{9 \lambda^2_{u}}}}{4 \lambda_{v}}\\
\approx  \lim_{\lambda_{u} < 0, \lambda_{v} \to 0} \frac{-3\lambda_{u} +3\lambda_{u}(1+\frac{1}{2}\frac{8 \lambda^2_{v}}{9 \lambda^2_{u}})}{4 \lambda_{v}} = \lim_{\lambda_{u} < 0, \lambda_{v} \to 0} \frac{\lambda_{v}}{3 \lambda_{u}} = 0.
\end{gathered}
\label{Eq:jump1} 
\end{align}
Analogously, if $\lambda_{u}(t) > 0$ for $ t\in [t_a,t_b]$, we have
\begin{align}
\begin{gathered}
\lim_{\lambda_{u} > 0, \lambda_{v} \to 0} \frac{-3\lambda_{u} - \sqrt{9 \lambda^2_{u} + 8 \lambda^2_{v}}}{4 \lambda_{v}} = \lim_{\lambda_{u} > 0, \lambda_{v} \to 0} \frac{-3\lambda_{u} - 3\lambda_{u} \sqrt{1+ \frac{8 \lambda^2_{v}}{9 \lambda^2_{u}}}}{4 \lambda_{v}}\\
\approx  \lim_{\lambda_{u} > 0, \lambda_{v} \to 0} \frac{-3\lambda_{u} - 3\lambda_{u}(1+\frac{1}{2}\frac{8 \lambda^2_{v}}{9 \lambda^2_{u}})}{4 \lambda_{v}} = \lim_{\lambda_{u} > 0, \lambda_{v} \to 0} -(\frac{3 \lambda_{u}}{2\lambda_{v}}+ \frac{\lambda_{v}}{3 \lambda_{u}}) = \lim_{\lambda_{u} > 0, \lambda_{v} \to 0} -\frac{3 \lambda_{u}}{2\lambda_{v}}.
\end{gathered}
\label{Eq:jump2} 
\end{align}

Therefore, from Eq.~(\ref{Eq:jump1}),
it is clear that the local minimum solution $\bar{\alpha}$ in Eq.~(\ref{Eq:optimal_alpha}) will automatically reduce to zero if $\lambda_{u}(t) < 0$ and $\lambda_{v}(t) = 0$ at some isolated points, indicating that Eq.~(\ref{Eq:optimal_alpha}) still holds in such case.
On the other hand, from Eq.~(\ref{Eq:jump2}), 
the sign of the local minimum solution $\bar{\alpha}$ will change  if $\lambda_{v}(t)$ crosses zero and $\lambda_{u}(t) > 0$ at the isolated points, which will result in the appearance of the discontinuous jump.
As a result, a discontinuous jump can be detected using the signs of $\lambda_{v}$ and $\lambda_{u}$.

Now we prove that the second case, that is, $\lambda_{v}(t) \equiv 0$ in a time interval $t \in [t_a,t_b]$, does not hold along an optimal trajectory. By contradiction, we assume that $\lambda_{v}(t) \equiv 0 $ for $t \in [t_a,t_b]$. Recall that $\lambda_\theta(t)\equiv 0$ along an optimal trajectory. Then, in view of Eq.~(\ref{EQ:costate_function}), we have
\begin{align}
\dot\lambda_{v}(t)=-2\frac{\lambda_{u}v(t)}{r(t)} \equiv 0, \forall~t\in [t_a,t_b].
\label{Eq:optimal_lambda}
\end{align}
Note that if $v(t)\equiv 0$ in any time interval, it will result in $\theta$ being constant in such time interval, as shown by Eq.~(\ref{EQ:dyna_equation}). This is obviously impossible during an orbital  transfer. Thus, to make Eq.~(\ref{Eq:optimal_lambda}) true, 
$\lambda_{u}(t)$ must be kept zero for $t \in [t_a,t_b]$. In this case, it implies
\begin{align}
\lambda_{r}(t)\equiv 0, \forall~t\in [t_a,t_b].
\label{Eq:optimal_r}
\end{align}
Clearly, if $\lambda_{v}(t) \equiv 0 $ for $t \in [t_a,t_b]$ does hold, the equation as follows will be valid
\begin{align}
\boldsymbol{\lambda}(t) = [\lambda_{r}(t), \lambda_{\theta}(t), \lambda_{u}(t),\lambda_{v}(t)] \equiv \boldsymbol{0}, \forall~t\in [t_a,t_b],
\label{Eq:optimal_00}
\end{align}
which contradicts the PMP. Thus, $\lambda_{v}(t)$ can only be zero at some isolated points.

\begin{figure}[!htp]
\centering
\includegraphics[width = 0.65\linewidth]{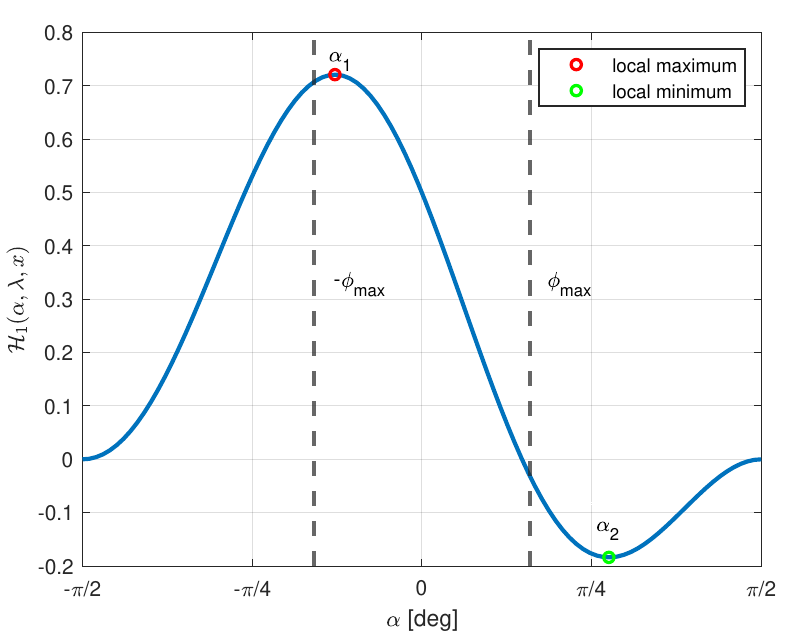}
\caption{A schematic view of $\mathscr H_1(\alpha, \boldsymbol{\lambda}, \boldsymbol{x})$ with $\lambda_{u}>0$ and $\lambda_{v}<0$.}
\label{Fig:polar_reference_app}
\end{figure}

Remember that the optimal guidance law is still ambiguous for the case that $\lambda_{v}(t)$ crosses zero and $\lambda_{u}(t) > 0$, as shown in Eq.~(\ref{Eq:jump2}). Because a nonlinear function takes its global minimum at one of its local minima or at one of the endpoints of its feasible domain \cite{oguri2022solar}, we have
\begin{align}
\alpha^* = \mathop {\rm{argmin}}\limits_{\alpha \in \{-\phi_{max},\bar{\alpha}, \phi_{max}\}} \mathscr H
\label{Eq:global_optimal_Ham}.
\end{align}
To further resolve the ambiguity in terms of $\alpha$, rewrite Eq.~(\ref{A1_H}) as
\begin{align}
\mathscr H(\alpha,\boldsymbol{\lambda},\boldsymbol{x}) = \mathscr H_1(\alpha,\boldsymbol{\lambda},\boldsymbol{x}) + \mathscr H_2(\boldsymbol{\lambda},\boldsymbol{x}),
\label{Eq:optimal_Ham}
\end{align}
where $\mathscr H_1$ is part of $\mathscr H$ related to $\alpha$, and $\mathscr H_2$ is the rest of $\mathscr H$ independent from $\alpha$. Then, $\mathscr H_1$ satisfies
\begin{align}
\mathscr H_1(\alpha,\boldsymbol{\lambda},\boldsymbol{x}) = \frac{\beta}{r^2}(\lambda_{u}\cos^3 \alpha + \lambda_{v} \sin \alpha \cos ^2 \alpha).
\label{Eq:optimal_Ham_1}
\end{align}
Recall Eq.~(\ref{EQ:att_constraints}) and $\phi_{max} \in (0,\frac{\pi}{2})$, thus we obtain
\begin{align}
\mathscr H_1(\frac{\pi}{2},\boldsymbol{\lambda},\boldsymbol{x}) = \mathscr H_1(-\frac{\pi}{2},\boldsymbol{\lambda},\boldsymbol{x}) = 0.
\label{Eq:optimal_Ham_2}
\end{align}
Substituting Eq.~(\ref{A1_alpha}) into Eq.~(\ref{Eq:optimal_Ham_1}) yields 
\begin{align}
\begin{cases}
   \mathscr H_{1}(\alpha_{1}) = \frac{\beta}{4r^{2}}\cos^{3}\alpha\left(  \lambda_{u}+\sqrt{9\lambda_{u}^{2}+8\lambda_{v}^{2}} \right)\geq 0,\\
   \mathscr H_{1}(\alpha_{2}) = \frac{\beta}{4r^{2}}\cos^{3}\alpha\left(  \lambda_{u} - \sqrt{9\lambda_{u}^{2}+8\lambda_{v}^{2}} \right)\leq 0,
\end{cases}
\end{align}
which indicates that $\alpha_1$ and $\alpha_2$ is the global maximum and global minimum for $\mathscr H$  in $[-\frac{\pi}{2}, \frac{\pi}{2}]$, respectively.
Define a variable $\Delta \mathscr H$ as
\begin{align}
\Delta \mathscr H =  \mathscr H_1(\phi_{max},\boldsymbol{\lambda},\boldsymbol{x}) - \mathscr H_1(-\phi_{max},\boldsymbol{\lambda},\boldsymbol{x}) =  \frac{2\beta}{r^2}\cos^2 \phi_{max} \sin \phi_{max} \lambda_{v}.
\label{Eq:optimal_Ham_3}
\end{align}
Without loss of generality, we consider the case that $\lambda_{v}<0$, which leads to a negative local maximum point $\alpha_1$ and a positive local minimum point $\alpha_2$. If $\alpha_2 < \phi_{max}$, then $\alpha_2$ is the global minimum point in $[-\phi_{max}, \phi_{max}]$. If $\alpha_2 > \phi_{max}$, then $\phi_{max}$ will be the global minimum in $[-\phi_{max}, \phi_{max}]$ since $\Delta \mathscr H < 0$, as shown in Fig.~\ref{Fig:polar_reference_app}. 
Hence, in view of Eq.~(\ref{Eq:optimal_alpha}), we have that the optimal guidance law $\alpha^*$ minimizing the Hamiltonian $\mathscr H$ w.r.t. $\alpha$ is 
\begin{align}
\alpha^*=\text{Median} {[-\phi_{max},\bar{\alpha},\phi_{max}]}, \text{where}~\bar{\alpha} = \arctan~\frac{-3 \lambda_{u}-\sqrt{9\lambda^{2}_{u}+8\lambda^{2}_{v}}}{4\lambda_{v}}.
\end{align}
   
Proof \textcolor{black}{of} Lemma \ref{LE:optimal_trajectory_lemma}. 
Taking Eq.~(\ref{para_initial}) into consideration, it is easy to see that the solution to the parameterized system in Eq.~(\ref{EQ:new_system}) at $\tau = 0$ represents a circular orbit with a radius of $R_0$. Moreover, regarding Eq.~(\ref{EQ:theta_law}), we have that Eq.~(\ref{EQ:tf_law}) holds at $\tau = 0$ for the pair $(\lambda_{R_0},\lambda_{U_0})$ and $\lambda_{V_0}$ determined by Eq.~(\ref{EQ:solve_equation_simple}). By propagating the parameterized system in Eq.~(\ref{EQ:new_system}) for $\tau \in [0,t_f]$, $\mathcal{F}$ meets all the necessary conditions for an optimal trajectory. By the definition of $\mathcal F_p$, it is obvious that $\mathcal F_p$ 
represents the solution space of an optimal trajectory, and $\tau \in [0,t_f]$ defines the time of flight. In other words, an optimal trajectory can be readily obtained simply by arbitrarily choosing two parameters $\lambda_{R_0}$ and $\lambda_{U_0}$ and solving an initial value problem governed by Eqs.~(\ref{EQ:new_system}) and (\ref{para_initial}). 



 
\bibliographystyle{elsarticle-num} 
\bibliography{myBib}
 




\end{document}